\documentclass[11pt,a4paper]{article}

\usepackage{amsmath,amssymb}
\newcommand{\e}{\varepsilon}
\newcommand{\va}{\varphi}
\newcommand{\D}{\Delta}
\newcommand{\La}{\Lambda}
\newcommand{\n}{\nabla}
\newcommand{\N}{\frac{N}{2}}
\newcommand{\NN}{\frac{N}{p}}

\newcommand{\p}{\partial}

\newcommand{\R}{\mathbb{R}}

\newcommand{\de}{\delta}
\newcommand{\del}{\bar{\delta}}
\newcommand{\ka}{\bar{\kappa}}

\newtheorem{definition}{Definition}
\newtheorem{theorem}{Theorem}
\newtheorem{proposition}{Proposition}

\newtheorem{remarka}{Remark}

\newtheorem{lemme}{Lemma}

\title{Cauchy problem for viscous shallow water equations  with a term of capillarity}
\author{Boris Haspot \thanks{Universit\'e Paris XII - Val de Marne 61, avenue
du G\'en\'eral de Gaulle 94 010 CRETEIL Cedex T\'el\'ephone : (33-1)
45 17 16 51 T\'el\'ecopie : (33-1) 45 17 16 49 e-mail :
haspot@univ-paris12.fr}}
\date{}
\begin{document}
\maketitle
\subsubsection*{Abstract}
In this article, we consider the compressible Navier-Stokes equation
with density dependent viscosity coefficients and a term of
capillarity introduced by Coquel et al in \cite{5CR}. This model
includes at the same time the barotropic Navier-Stokes equations with variable viscosity coefficients, shallow-water system
and the model of Rohde in \cite{5Ro}.\\
We first study the well-posedness of the model in 
critical regularity spaces with respect to the scaling of the
associated equations. In a functional setting as close as possible
to the physical energy spaces, we prove global existence of
solutions close to a stable equilibrium, and local in time existence
for solutions with general initial data. Uniqueness is also
obtained.
\section{Introduction}
This paper is devoted to the Cauchy problem for the compressible
Navier-Stokes equation with viscosity coefficients depending on the
density and with a capillary term coming from the works of Coquel,
Rohde and theirs collaborators in \cite{5CR}, \cite{5Ro}. Let $\rho$
and $u$ denote the density and the velocity of a compressible
viscous fluid. As usual, $\rho$ is a non-negative function and $u$
is a vector valued function defined on $\R^{N}$. Then, the
Navier-Stokes equation for compressible fluids endowed with internal
capillarity introduced in \cite{5Ro} reads:
$$
\begin{cases}
\begin{aligned}
&\p_{t}\rho+{\rm div}(\rho u)=0,\\
&\p_{t}(\rho u)+{\rm div}(\rho u\otimes u)-{\rm
div}(2\mu(\rho)Du)-\n(\lambda(\rho){\rm div}u)
+\n P(\rho)=\kappa\rho\n D[\rho],\\
\end{aligned}
\end{cases}
\leqno{(SW)}
$$
supplemented by the initial  condition:
$$\rho_{/ t=0}=\rho_{0},\;\;\rho u_{/ t=0}=\rho_{0}u_{0}$$
and:
$$D[\rho]=\phi*\rho-\rho$$
where $\phi$ is chosen so that:
$$\phi\in L^{\infty}(\R^{N})\cap C^{1}(\R^{N})\cap W^{1,\,1}(\R^{N}),\;\;\;\int_{\R^{N}}\phi(x)dx=1,\;\;\phi
\;\;\mbox{even},\;\mbox{and}\;\;\phi\geq0$$ and where $P(\rho)$
denotes the pressure, $\mu$ and $\lambda$ are the two Lam\'e
viscosity coefficients (they depend regularly on the density $\rho$)
satisfying:
$$\mu>0\;\;\;\;\;\;2\mu+N\lambda\geq0.$$
($\mu$ is sometimes called the shear viscosity of the fluid, while
$\lambda$ is usually referred to as the second
viscosity coefficient).\\
several physical models arise as a particular case of system $(SW)$:
\begin{itemize}
\item when $\kappa=0$ $(SW)$
represents  compressible Navier-Stokes model with variable viscosity
coefficients.
\item when $\kappa=0$ and $\mu(\rho)=\rho$, $\lambda(\rho)=0$,
$P(\rho)=\rho^{2}$, $N=2$ then $(SW)$ describes the system of
shallow-water.
\item when $\kappa=0$ and $\mu$, $\lambda$ are constant, $(SW)$
reduce to the Rohde model of chapter four.
\end{itemize}
One of the major difficulty of compressible fluid mechanics is to
deal with the vacuum. The problem of existence of global solution in
time for Navier-Stokes equations was addressed in one dimension for
smooth enough data by Kazhikov and Shelukin in \cite{5K}, and for
discontinuous ones, but still with densities away from zero, by Serre
in  \cite{5S} and Hoff in \cite{5H1}. Those results have been
generalized to higher dimension by Matsumura and Nishida in
\cite{5MN} for smooth data close
to equilibrium and by Hoff in the case of discontinuous data in \cite{5H2,5H3}.\\
Concerning large initial data, Lions showed in \cite{5L2} the global
existence of weak solutions for $\gamma\geq\frac{3}{2}$ for $N=2$
and $\gamma\geq\frac{9}{5}$ for $N=3$. Let us mention that Feireisl
in \cite{5F}  generalized the result to $\gamma>\frac{N}{2}$ in
establishing that we can obtain renormalized solution without
imposing that $\rho\in L^{2}_{loc}$, for this he introduces the
concept of oscillation defect measure evaluating the
loss of compactness.\\
Other results provide the full range $\gamma>1$ under symmetries
assumptions on the initial datum, see Jiang and Zhang \cite{5J}. All
those results do not require to be far from the vacuum. However they
rely strongly on the assumption that the
viscosity coefficients are bounded below by a positive constant.
This non
physical assumption allows to get some estimates on the gradient of the velocity field.\\
The main difficulty when dealing with vanishing viscosity
coefficients on vacuum is that the velocity cannot even be defined
when the density vanishes and so we cannot use some properties of
parabolicity of the momentum equation, see \cite{5CK1},
\cite{5CK2}.\\
The first result handling this difficulty is due to Bresch,
Desjardins and Lin in \cite{5BDL}. They show the existence of global
weak solution for Korteweg system in choosing specific type of
viscosity where $\mu$ and $\lambda$ are linked.\\
The result was later improved by Bresch and Desjardins in \cite{5BD}
to include the case of vanishing capillarity ($\kappa=0$), but with
an additional quadratic friction term $r\rho u|u|$ (see also
\cite{5BD2}). However, those estimates are not enough to treat the
case without capillarity and friction effects $\kappa=0$ and $r=0$
(which corresponds to equation
(1) with $h(\rho)=\rho$ and $g(\rho)=0$).\\
The main difficulty, to prove the stability of $(SW)$ , is to pass
to the limit in the term $\rho u\otimes u$ (which requires the
strong convergence of $\sqrt{\rho} u$). Note that this is easy when
the viscosity coefficients are bounded below by a positive constant.
On the other hand, the new bounds on the gradient of the density
make the control of the pressure term far simpler
than in the case of constant viscosity coefficients.\\
In \cite{5BD2} Bresch and Desjardins show a result of global
existence of weak solution for the non isothermal Navier-Stokes
equation by imposing some condition between the viscosity
coefficient and  a bound by below on the viscosity
coefficient. A. Mellet and A. Vasseur in using  the same entropy
inequality get a very interesting result of stability. They get more
general estimate, which hold for any viscosity coefficients
$\mu(\rho)$, $\lambda(\rho)$ satisfying the relation:
\begin{equation}
\mu(\rho)=\rho\lambda^{'}(\rho)-\lambda(\rho). \label{5coeff}
\end{equation}
Mellet and Vasseur show in \cite{5MV} the $L^{1}$ stability of weak
solutions of Navier-Stokes compressible isotherm
 under some conditions on the
viscosity coefficients (including (\ref{5coeff})) but without any
additional regularizing terms. The interest of this result is to
consider conditions where the viscosity coefficients vanish on the
vacuum set. It includes the case $\mu(\rho)=\rho$, $\lambda(\rho)=0$
(when $N=2$ and $\gamma=2$, where we recover the Saint-Venant model
for Shallow water). The key to the proof is  a new energy inequality
on the velocity and a gain of integrability, which
allows to pass to the limit.\\
The existence and uniqueness of local classical solutions for $(SW)$
with smooth initial data such that the density $\rho_{0}$ is bounded
and bounded away from zero (i.e.,
$0<\underline{\rho}\leq\rho_{0}\leq M$)
has been stated by Nash in \cite{5Na}. Let us emphasize that no stability condition was required there.\\
On the other hand, for small smooth perturbations of a stable
equilibrium with constant positive density, global well-posedness
has been proved in \cite{5MN}. Many recent works have been devoted
to the qualitative behavior of solutions for large time (see for
example \cite{5H1,5K}). Refined functional analysis has been used
for the last decades, ranging from Sobolev, Besov, Lorentz and
Triebel spaces to describe the regularity and long time behavior of
solutions to the compressible model \cite{5So}, \cite{5V},
\cite{5H4}, \cite{5K1}. The most important result on the system of
Navier-Stokes compressible isothermal comes from R. Danchin in
\cite{5DG} and \cite{5DL} who show the existence of global solution
and uniqueness with initial data close from a equilibrium, and he
has the same result in finite time. The interest is that he can work
in {\it critical} Besov space ({\it critical} in the sense
of the scaling of the equation.)
\\
We generalize the result of R. Danchin in considering general
viscosity coefficient and in connecting this result with those of A.
Mellet and A. Vasseur. This result improves too the case of strong
solution for the shallow-water system, where W.Wang and C-J Xu in
\cite{5W} have got global existence in time for small initial
data with $h_{0},\;u_{0}\in H^{2+s}$ with $s>0$.
\subsection{Notations and main results}
We will mainly consider the global well-posedness problem for
initial data close enough to stable equilibria. Here we want to
investigate the well-posedness of the system $(SW)$ problem in
critical spaces, that is, in spaces which are invariant by the
scaling of the equations.
Let us explain precisely the scaling of the system. We can easily
check that, if $(\rho,u)$ solves $(SW)$, so does
$(\rho_{\lambda},u_{\lambda})$, where:
$$\rho_{\lambda}(t,x)=\rho(\lambda^{2}t,\lambda x)\,\,\,\,\mbox{and}\,\,\,\,u_{\lambda}(t,x)=
\lambda u(\lambda^{2}t,\lambda x)$$
\\
provided the pressure law $P$ has been changed into
$\lambda^{2}P$.
\begin{definition}
We say that a functional space is critical with respect to the
scaling of the equation if the associated norm is invariant under
the transformation:
$$(\rho,u)\longrightarrow(\rho_{\lambda},u_{\lambda})$$
(up to a constant independent of $\lambda$).
\end{definition}
This suggests us to choose initial data $(\rho_{0},u_{0})$ in spaces
whose norm is invariant for all $\lambda>0$ by
$(\rho_{0},u_{0})\longrightarrow(\rho_{0}(\lambda\cdot),\lambda
u_{0}(\lambda\cdot)).$
\\
\\
A natural candidate is the homogeneous Sobolev space
$\dot{H}^{N/2}\times (\dot{H}^{N/2-1})^{N}$, but since
$\dot{H}^{N/2}$ is not included in $L^{\infty}$, we cannot expect to
get $L^{\infty}$ control on the density when
$\rho_{0}\in\dot{H}^{N/2}$.
\\
This is the reason why as in the chapter two, instead of the
classical homogeneous Sobolev space, we will consider homogeneous
Besov spaces $B^{N/2}_{2,1}\times (B^{N/2-1}_{2,1})^{N}$ with the
same derivative index.
This allows to control the density from below and from above,
without requiring more regularity on derivatives of $\rho$.
In the sequel, we will need to control the vacuum, this motivates
the following definition:
\begin{definition}
\label{5bar} Let $\bar{\rho}>0$, $\bar{\theta}>0$. We will note in
the sequel:
$$q=\frac{\rho-\bar{\rho}}{\bar{\rho}}$$
\end{definition}
Let us first state a result of global existence and uniqueness of
$(SW)$ for initial data close to a equilibrium.
\begin{theorem} \label{5theo1}
Let $N\geq2$. Let $\bar{\rho}>0$ be such that:
$P^{'}(\bar{\rho})>0$, $\mu(\bar{\rho})>0$ and $2\mu(\bar{\rho})+\lambda(\bar{\rho})>0.$
There exist two positive constants $\e_{0}$ and $M$ such that if
$q_{0}\in\widetilde{B}^{\frac{N}{2}-1,\N}$, $u_{0}\in B^{\N-1}$ and:
$$\|q_{0}\|_{\widetilde{B}^{\frac{N}{2}-1,\N}}+\|u_{0}\|_{B^{\N-1}}\leq\e_{0}$$
then $(SW)$ has a unique global solution $(q,u)$ in $E^{\N}$ which
satisfies:
$$\|(q,u)\|_{E^{\N}}\leq M \big(\,\|q_{0}\|_{\widetilde{B}^{\frac{N}{2}-1,\N}}+\|u_{0}\|_{B^{\N-1}}\big),$$
for some $M$ independent of the initial data where:
$$
\begin{aligned}
E^{\N}=[C_{b}(\R^{+},\widetilde{B}^{\N-1,\N})\cap L^{1}(\R^{+},\widetilde{B}^{\N+1,\N})]\times[C_{b}(\R_{+}&,B^{\N-1})^{N}\\
&\cap L^{1}(\R^{+},B^{\N+1})^{N}].\\
\end{aligned}
$$
\end{theorem}
In the following theorems, we want to show some result of existence
and uniqueness in finite time for large initial velocity and
initial density close to some constant.
\\
The following result shows the existence and uniqueness in finite
time with initial data in critical Besov space for the scaling of
the equations. However as we said, we need for some technical
reasons of an hypothesis of smallness on $q_{0}$. We note that we
work on Besov space $B^{s}_{p}$ with general index $p$ on the
integrability.
\begin{theorem}
\label{5theo2} Let $p\in[1,+\infty[$.
Let $q_{0}\in B^{\NN}_{p}$ and $u_{0}\in B^{\NN-1}_{p}$. Assume also
that: $\|q_{0}\|_{B^{\NN}_{p}}\leq\e$
for a suitably small positive constant $\e>0$.\\
Under the assumptions of the theorem \ref{5theo1} for the physical coefficients, there exists a time $T>0$ such that the following results hold:
\begin{enumerate}
\item Existence: If $p\in[1,2N[$ then system $(SW)$ has a solution $(q,u)$ in $F^{\NN}_{p}$ with:
$$F^{\NN}_{p}=\widetilde{C}_{T}(B^{\NN}_{p})\times\big(L^{1}_{T}(B^{\NN+1}_{p})
\cap\widetilde{C}_{T}(B^{\NN-1}_{p}\big).$$
\item Uniqueness: If in addition $1\leq p\leq N$ then uniqueness holds in $F^{\NN}_{p}$.
\end{enumerate}
Moreover we have a control on the time $T$ which may be bounded from
below by:
$$\min\biggl(\eta,\max\big(t>0,\sum_{q\in\mathbb{Z}}2^{q(\NN-1)}\|\D_{q}u_{0}\|_{L^{p}}(\frac{1
-e^{-c\widetilde{\nu}te^{2q}}}{c\widetilde{\nu}})\leq\frac{\e\widetilde{\nu}^{2}}{\widetilde{\nu}+U_{0}}\big)\biggl).$$
where
$\widetilde{\nu}=\min(\mu(\bar{\rho}),\lambda(\bar{\rho})+2\mu(\bar{\rho}))$
and $U_{0}=\|u_{0}\|_{B^{\NN}}$.
\end{theorem}
In the next theorem, we consider the case when the initial
variational density $q_{0}$ belongs to
$\bar{\rho}+\widetilde{B}^{\NN,\NN+\e}_{p}$ with $\e>0$ and
satisfies $0<\underline{\rho}<\rho$. Here we do not suppose a
smallness condition on $\|q_{0}\| _{B^{\NN}}$ but we impose more
regularity.
\begin{theorem}
\label{5theo3} Let $\e\in]0,1[$ and $p\in[1,\frac{N}{1-\e}[$ and we assume that the physical coefficients verify
the same hypothesis as in theorem \ref{5theo1}.
Let $\rho_{0}\in\bar{\rho}+\widetilde{B}^{\NN,\NN+\e}_{p}$ for a
constant $\bar{\rho}>0$, $u_{0}\in B^{\NN+\e-1}$. Assume that there
is a constant $\underline{\rho}$ such that:
$$0<\underline{\rho}\leq\rho_{0}.$$
There exists a time $T>0$ such that system $(SW)$ has a unique solution $(q,u)$ in $F^{\NN}_{p+\e}$.
Moreover we have a control on the time $T$ which may be bounded from
below by:
$$\min\biggl(\eta,\max\big(t>0,\sum_{q\in\mathbb{Z}}2^{q(\NN-1)}\|\D_{q}u_{0}\|_{L^{p}}(\frac{1
-e^{-c\widetilde{\nu}te^{2q}}}{c\widetilde{\nu}})\leq\frac{\e\widetilde{\nu}^{2}}{\widetilde{\nu}+U_{0}}\big)\biggl).$$
where:
$\widetilde{\nu}=\min(\mu(\bar{\rho}),\lambda(\bar{\rho})+2\mu(\bar{\rho}))$
and $U_{0}=\|u_{0}\|_{B^{\NN}}$.
\end{theorem}
The present chapter is structured as follows.\\
In the Section \ref{5section3}, we recall some basic facts about
Littlewood-Paley
decomposition and Besov spaces.\\
In the Section \ref{5section4} we prove the theorem \ref{5theo1} and
in the Section \ref{5section5} we show the theorem \ref{5theo2}. We
conclude in the Section \ref{5section6} by the proof of the theorem
\ref{5theo3}. 
\section{Littlewood-Paley theory and Besov spaces}
\label{5section3}
\subsection{Littlewood-Paley decomposition}
Littlewood-Paley decomposition  corresponds to a dyadic
decomposition  of the space in Fourier variables.\\
We can use for instance any $\varphi\in C^{\infty}(\R^{N})$,
supported in
${\cal{C}}=\{\xi\in\R^{N}/\frac{3}{4}\leq|\xi|\leq\frac{8}{3}\}$
such that:
$$\sum_{l\in\mathbb{Z}}\varphi(2^{-l}\xi)=1\,\,\,\,\mbox{if}\,\,\,\,\xi\ne 0.$$
Denoting $h={\cal{F}}^{-1}\varphi$, we then define the dyadic
blocks by:
$$\D_{l}u=\varphi(2^{-l}D)u=2^{lN}\int_{\R^{N}}h(2^{l}y)u(x-y)dy\,\,\,\,\mbox{and}\,\,\,S_{l}u=\sum_{k\leq
l-1}\D_{k}u\,.$$ Formally, one can write that:
$$u=\sum_{k\in\mathbb{Z}}\D_{k}u\,.$$
This decomposition is called homogeneous Littlewood-Paley
decomposition. Let us observe that the above formal equality does
not hold in ${\cal{S}}^{'}(\R^{N})$ for two reasons:
\begin{enumerate}
\item The right hand-side does not necessarily converge in
${\cal{S}}^{'}(\R^{N})$.
\item Even if it does, the equality is not
always true in ${\cal{S}}^{'}(\R^{N})$ (consider the case of the
polynomials).
\end{enumerate}
However, this equality holds true modulo polynomials hence
homogeneous Besov spaces will be defined modulo the polynomials,
according to \cite{5DG}.
\subsection{Homogeneous Besov spaces and first properties}
\begin{definition}
For
$s\in\R,\,\,p\in[1,+\infty],\,\,q\in[1,+\infty],\,\,\mbox{and}\,\,u\in{\cal{S}}^{'}(\R^{N})$
we set:
$$\|u\|_{B^{s}_{p,q}}=(\sum_{l\in\mathbb{Z}}(2^{ls}\|\D_{l}u\|_{L^{p}})^{q})^{\frac{1}{q}}.$$
\end{definition}
A difficulty due to the choice of homogeneous spaces arises at this
point. Indeed, $ \|.\|_{B^{s}_{p,q}}$ cannot be a norm on
$\{u\in{\cal{S}}^{'}(\R^{N}),\|u\|_{B^{s}_{p,q}}<+\infty\}$ because
$\|u\|_{B^{s}_{p,q}}=0$ means that $u$ is a polynomial. This
enforces us to adopt the following definition for homogeneous Besov
spaces, see \cite{5DG}.
\begin{definition}
Let $s\in\R,\,\,p\in[1,+\infty],\,\,q\in[1,+\infty]$.\\
Denote $m=[s-\NN]$ if $s-\NN\notin\mathbb{Z}$ or $q>1$ and
$m=s-\NN-1$ otherwise.
\begin{itemize}
\item If $m<0$, then we define $B^{s}_{p,q}$ as:
$$B^{s}_{p,q}=\biggl\{u\in{\cal{S}}^{'}(\R^{N})\big/\;\|u\|_{B^{s}_{p,q}}<\infty\,\,\,\mbox{and}
\,\,\,u=\sum_{l\in\mathbb{Z}}\D_{l}u\,\,\mbox{in}\,\,{\cal{S}}^{'}(\R^{N})\biggl\}\,.$$
\item If $m\geq 0$, we denote by ${\cal{P}}_{m}[\R^{N}]$ the set of
polynomials of degree less than or equal to $m$ and we set:
$$B^{s}_{p,q}=\biggl\{u\in{\cal{S}}^{'}(\R^{N})/{\cal{P}}_{m}[\R^{N}]\;\big/\;
\|u\|_{B^{s}_{p,q}}<\infty\,\,\,and\,\,\,u=\sum_{l\in\mathbb{Z}}
\D_{l}u\,\,in\,\,{\cal{S}}^{'}(\R^{N}){\cal{P}}_{m}[\R^{N}]\biggl\}\,.$$
\end{itemize}
\end{definition}
The definition of $B^{s}_{p,r}$ does not depend on the choice of the Littlewood-Paley decomposition.
\begin{remarka}
In the sequel, we will use only Besov space $B^{s}_{p,q}$ with $q=1$
and we will denote them by $B^{s}_{p}$ or even by $B^{s}$ if there
is no ambiguity on the index $p$.\\
Let us now state some basic properties for those Besov spaces.
\end{remarka}
\begin{proposition}
\label{derivation} The following properties holds:
\begin{enumerate}
\item Density: If $p<+\infty$ and $|s|\leq N/p$ , then
$C^{\infty}_{0}$ is dense in $B^{s}_{p}$.
\item Derivatives: there exists a constant universal $C$
such that:
$$C^{-1}\|u\|_{B^{s}_{p,r}}\leq\|\n u\|_{B^{s-1}_{p,r}}\leq
C\|u\|_{B^{s}_{p,r}}.$$
\item Sobolev embeddings: If
$p_{1}<p_{2}$ and $r_{1}\leq r_{2}$ then $B^{s}_{p_{1},r_{1}}\hookrightarrow
B^{s-N(\frac{1}{p_{1}}-\frac{1}{p_{2}})}_{p_{2},r_{2}}$.
\item Algebraic
properties: For $s>0$, $B^{s}_{p,r}\cap L^{\infty}$ is an algebra. Moreover, for any $p\in[1,+\infty]$ then
$B^{\NN}_{p,1}\hookrightarrow B^{\NN}_{p,\infty}\cap L^{\infty}$, and $B^{\NN}_{p,1}$ is an algebra if $p$ is finite.
\item Real interpolation:
 $(B^{s_{1}}_{p,r},B^{s_{2}}_{p,r})_{\theta,r^{'}}=B^{\theta
s_{1}+(1-\theta)s_{2}}_{p,r^{'}}$.
\end{enumerate}
\end{proposition}
\subsection{Hybrid Besov spaces and Chemin-Lerner spaces}
Hybrid  Besov spaces are functional spaces where regularity
assumptions are different in low frequency and high frequency, see
\cite{5DG}.
We are going to give the definition of this this new spaces and give
some of their main properties.
\begin{definition}
Let $s,\,t\in\R$.We set:
$$\|u\|_{\widetilde{B}^{s,t}_{p,r}}=\big(\sum_{q\leq 0}(2^{qs}\|\D_{q}u\|_{L^{p}})^{r}\big)^{\frac{1}{r}}+\big(\sum_{q>
0}(2^{qt}\|\D_{q}u\|_{L^{p}})^{r}\big)^{\frac{1}{r}}\,.$$
Denote $m=[s-\NN]$ if $s-\NN\notin\mathbb{Z}$ or $r>1$ and
$m=s-\NN-1$ otherwise, we then define:
\begin{itemize}
\item $\widetilde{B}^{s,t}_{p}=\{u\in{\cal{S}}^{'}(\R^{N})\;\big/\;\|u\|_{\widetilde{B}^{s,t}_{p}}<+\infty\}$,
if $m<0$
\item $\widetilde{B}^{s,t}_{p}=\{u\in{\cal{S}}^{'}(\R^{N})/{\cal{P}}_{m}[\R^{N}]\;\big/\;\|u\|_{\widetilde{B}^{s,t}_{p}}<+\infty\}$
if $m\geq 0$.
\end{itemize}
\end{definition}
Let now give some properties of these hybrid spaces and some results
on how they behave with respect to the product. The following
results come directly from the paradifferential calculus.
\begin{proposition}
\label{important} We give here some results of inclusion:
\begin{enumerate}
\item We have $\widetilde{B}^{s,s}_{p,r}=B^{s}_{p,r}$.
\item If $s\leq t$ then $\widetilde{B}^{s,t}_{p,r}=B^{s}_{p,r}\cap B^{t}_{p,r}$  or
if $s>t$ then $\widetilde{B}^{s,t}_{p,r}=B^{s}_{p,r}+B^{t}_{p,r}$.
\item If $s_{1}\leq s_{2}$ and $t_{1}\geq t_{2}$ then
$\widetilde{B}^{s_{1},t_{1}}_{p,r}\hookrightarrow
\widetilde{B}^{s_{2},t_{2}}_{p,r}$.
\end{enumerate}
\end{proposition}
\begin{proposition}
\label{5produit1}
For all $s,\,t>0$, $1\leq r,p\leq +\infty$, the following inequality holds true:
\begin{equation}
\|uv\|_{\widetilde{B}^{s,t}_{p,r}}\leq
C(\|u\|_{L^{\infty}}\|v\|_{\widetilde{B}^{s,t}_{p,r}}+\|v\|_{L^{\infty}}\|u\|_{\widetilde{B}^{s,t}_{p,r}})\,.
\label{5prodinteressant1}
\end{equation}
For all $s_{1},\,s_{2},\,t_{1},\,t_{2}\leq \frac{N}{p}$ such that
$\min(s_{1}+s_{2},\,t_{1}+t_{2})>0$ we have:
\begin{equation}
\|uv\|_{\widetilde{B}^{s_{1}+t_{1}-\NN,s_{2}+t_{2}-\NN}_{p,r}}\leq
C\|u\|_{\widetilde{B}^{s_{1},t_{1}}_{p,r}}\|v\|_{\widetilde{B}^{s_{2},t_{2}}_{p,\infty}}\,.
\label{5prodinteressant2}
\end{equation}
\begin{equation}
\|uv\|_{B^{s}_{p,r}}\leq C\|u\|_{B^{s}_{p,r}}\|v\|_{B^{\NN}_{p,\infty}\cap L^{\infty}}\;\;\;\mbox{if}\;|s|<\NN.
\label{5prodinteressant3}
\end{equation}
\end{proposition}
For a proof of this proposition see \cite{5DG}. The limit case $s_{1}+s_{2}=t_{1}+t_{2}=0$ in (\ref{5prodinteressant2}) is of interest.
When $p\geq2$, the following estimate holds true whenever $s$ is in the range $(-\NN,\NN]$ (see e.g. \cite{5RS}):
\begin{equation}
\|uv\|_{B^{-\NN}_{p,\infty}}\leq C\|u\|_{B^{s}_{p,1}}\|v\|_{B^{-s}_{p,\infty}}.
\label{5limite}
\end{equation}
The study of non stationary PDE's requires space of type $L^{\rho}(0,T,X)$ for appropriate Banach spaces $X$. In our case, we
expect $X$ to be a Besov space, so that it is natural to localize the equation through Littlewood-Payley decomposition. But, in doing so, we obtain
bounds in spaces which are not type $L^{\rho}(0,T,X)$ (except if $r=p$).
We are now going to
define the spaces of Chemin-Lerner in which we will work (see \cite{5CT}), which are
a refinement of the spaces
$L_{T}^{\rho}(B^{s}_{p,r})$.
$\hspace{15cm}$
\begin{definition}
Let $\rho\in[1,+\infty]$, $T\in[1,+\infty]$ and $s_{1},s_{2}\in\R$. We then
denote:
$$\|u\|_{\widetilde{L}^{\rho}_{T}(\widetilde{B}^{s_{1},s_{2}}_{p,r})}=
\big(\sum_{l\leq0}2^{lrs_{1}}(\|\D_{l}u(t)\|_{L^{\rho}_{T}(L^{p})}^{r}\big)^{\frac{1}{r}}+
\big(\sum_{l>0}2^{lrs_{2}}(\int_{0}^{T}\|\D_{l}u(t)\|^{\rho}_{L^{p}}dt)^{\frac{r}{\rho}})\big)^{\frac{1}{r}}\,.$$
\end{definition}
We note that thanks to Minkowsky inequality we have:
$$
\begin{aligned}
&\|u\|_{L^{\rho}_{T}(\widetilde{B}^{s_{1},s_{2}}_{p,r})}\leq\|u\|_{\widetilde{L}^{\rho}_{T}
(\widetilde{B}^{s_{1},s_{2}}_{p,r})}\;\;\;\mbox{if}\;\;\rho\leq r,\\
&\|u\|_{\widetilde{L}^{\rho}_{T}
(\widetilde{B}^{s_{1},s_{2}}_{p,r})}\leq\|u\|_{L^{\rho}_{T}(\widetilde{B}^{s_{1},s_{2}}_{p,r})}
\;\;\;\mbox{if}\;\;\rho\geq r.\\
\end{aligned}
$$
We then define the space:
$$\widetilde{L}^{\rho}_{T}(\widetilde{B}^{s_{1},s_{2}}_{p})=\{u\in
L^{\rho}_{T}(\widetilde{B}^{s_{1},s_{2}}_{p})/\|u\|_{\widetilde{L}^{\rho}_{T}(\widetilde{B}^{s_{1},s_{2}}_{p})}<\infty\}\,.$$
We
denote moreover by $\widetilde{C}_{T}(\widetilde{B}^{s_{1},s_{2}}_{p})$
the set of those functions of
$\widetilde{L}^{\infty}_{T}(\widetilde{B}^{s_{1},s_{2}}_{p})$ which are
continuous from $[0,T]$ to $\widetilde{B}^{s_{1},s_{2}}_{p}$.
In the
sequel we are going to give some properties of this spaces
concerning the interpolation and their relationship with the heat
equation.
\begin{proposition}
\label{5interp} Let $s,\,t,\,s_{1},\,s_{2}\in\R$,
$r,\,\rho,\,\rho_{1},\,\rho_{2}\in[1,+\infty]$. We have:
\begin{enumerate}
\item Interpolation:
$$\|u\|_{\tilde{L}^{\rho}_{T}(\tilde{B}^{s,t}_{p,r})}\leq\|u\|_{\tilde{L}
^{\rho_{1}}_{T}(\tilde{B}^{s_{1},t_{1}}_{p,r})}^{\theta}\|u\|_{\tilde{L}^
{\rho_{2}}_{T}(\tilde{B}^{s_{2},t_{2}}_{p,r})}^{1-\theta}$$ with
$\frac{1}{\rho}=\frac{\theta}{\rho_{1}}+\frac{1-\theta}{\rho_{2}}$,
$s=\theta s_{1}+(1-\theta)s_{2}$ and $t=\theta
t_{1}+(1-\theta)t_{2}$.
\item Embedding:
$$\tilde{L}^{\rho}_{T}(\tilde{B}^{s,t}_{p})\hookrightarrow
L^{\rho}_{T}(C_{0})\;\;\mbox{and}\;\;\tilde{C}_{T}(B^{\NN}_{p})\hookrightarrow
C([0,T]\times\R^{d})$$
\end{enumerate}
\end{proposition}
Here we recall a result of interpolation which explains the link
of the space $B^{s}_{p,1}$ with the homogeneous spaces , see
\cite{5DFourier}.
\begin{proposition}
\label{5interp1}
There exists a constant $C$ such that for all $s\in\R$, $\e>0$ and
$1\leq p\leq+\infty$, we have
$$\|u\|_{B^{s}_{p,1}}\leq C\frac{1+\e}{\e}\|u\|_{B^{s}_{p,\infty}}\biggl(1+\log\frac{\|u\|_{B^{s-\e}_{p,\infty}}+
\|u\|_{B^{s+\e}_{p,\infty}}}
{\|u\|_{B^{s}_{p,\infty}}}\biggl).$$ \label{5Yudov}
\end{proposition}
To finish with we adapt the results of the paradifferential calculus
on the product of Besov function to the spaces of Chemin-Lerner. So
we have the following properties:
\begin{proposition}
\label{5produit}
Let $p,r\in[1,+\infty]$. We have the two following properties:
\begin{itemize}
\item Let $s>0,\;t>0$,
$1/\rho_{2}+1/\rho_{3}=1/\rho_{1}+1/\rho_{4}=1/\rho\leq1$,
$u\in\tilde{L}^{\rho_{3}}_{T}(\tilde{B}^{s,t}_{p,r})\cap
\tilde{L}^{\rho_{1}}_{T}(L^{\infty})$ and
$v\in\tilde{L}^{\rho_{4}}_{T}(\tilde{B}^{s,t}_{p,r})\cap
\tilde{L}^{\rho_{2}}_{T}(L^{\infty})$. Then
$uv\in\tilde{L}^{\rho}_{T}(\tilde{B}^{s,t}_{p,r})$ and we have:
$$\|uv\|_{\tilde{L}^{\rho}_{T}(\tilde{B}^{s,t}_{p,r})}\lesssim
\|u\|_{\tilde{L}^{\rho_{1}}_{T}(L^{\infty})}\|v\|_{\tilde{L}^{\rho_{4}}_{T}
(\tilde{B}^{s,t}_{p,r})}+\|v\|_{\tilde{L}^{\rho_{2}}_{T}(L^{\infty})}
\|u\|_{\tilde{L}^{\rho_{3}}_{T}(\tilde{B}^{s,t}_{p,r})}$$
\item If $s_{1},s_{2},t_{1},t_{2}\leq \NN$ , $s_{1}+s_{2}>0$,
$t_{1}+t_{2}>0$ $1/\rho_{1}+1/\rho_{2}=1/\rho\leq1$,
$u\in\tilde{L}^{\rho_{1}}_{T}(B^{s_{1},t_{1}}_{p,r})$ and
$v\in\tilde{L}^{\rho_{2}}_{T}(B^{s_{2},t_{2}}_{p,r})$ then
$uv\in\tilde{L}^{\rho}_{T}(B^{s_{1}+s_{2}-d/2}_{2})$ and
$$\|uv\|_{\tilde{L}^{\rho}_{T}(B^{s_{1}+s_{2}-\NN,t_{1}+t_{2}-\NN}_{p,r})}\lesssim\|u\|_{\tilde{L}^{\rho_{1}}_{T}
(B^{s_{1},t_{1}}_{p,r})}\|v\|_{\tilde{L}^{\rho_{2}}_{T}(B^{s_{2},t_{2}}_{p,r})}\,.$$
\end{itemize}
\end{proposition}
The analogous of the endpoint estimate (\ref{5limite}) reads (for $p\geq2$):
\begin{equation}
\|uv\|_{\widetilde{L}^{\rho}_{T}(B^{-\NN}_{p,\infty})}\lesssim\|u\|_{\widetilde{L}^{\rho_{1}}_{T}(B^{s}_{p,1})}
\|u\|_{\widetilde{L}^{\rho_{2}}_{T}(B^{-s}_{p,\infty})},
\end{equation}
whenever $s$ is in the range $(-\NN,NN]$ and $\frac{1}{\rho_{1}}+\frac{1}{\rho_{2}}=\frac{1}{\rho}\leq1$ (see the proof in \cite{5DD}).
For a proof of this proposition see \cite{5DG}. Finally we need an
estimate on the composition of functions in the spaces
$\tilde{L}^{\rho}_{T}(\tilde{B}^{s}_{p})$.
\begin{proposition}
\label{5composition} Let $s>0$, $r\in[1,+\infty]$ and $F\in W^{s+2,\infty}_{loc}(\R^{N})$ such that $F(0)=0$.
There exists a function $C$ depending only on $s$, $p$, $N$ and $F$, and such that:
$$\|F(u)\|_{\tilde{L}^{\rho}_{T}(\tilde{B}^{s_{1},s_{2}}_{p,r})}\leq C
(\|u\|_{L^{\infty}_{T}(L^{\infty})})\|u\|_{\tilde{L}^{\rho}_{T}(\tilde{B}^{s_{1},s_{2}}_{p,r})}.$$
If $v,\,u\in\tilde{L}^{\rho}_{T}(B^{s_{1},s_{2}}_{p})\cap
L^{\infty}_{T}(L^{\infty})$ and $G\in
W^{[s]+3,\infty}_{loc}(\R^{N})$ then $G(u)-G(v)$ belongs to
$\tilde{L}^{\rho}_{T}(B^{s_{1},s_{2}}_{p})$ and it exists a constant
C depending only of $s, p ,N\;\mbox{and}\;G$ such that:
$$
\begin{aligned}
&\|G(u)-G(v)\|_{\tilde{L}^{\rho}_{T}(B^{s_{1},s_{2}}_{p,r})}\leq
\,\,C(\|u\|_{L^{\infty}_{T}(L^{\infty})},\|v\|_{L^{\infty}_{T}(L^{\infty})})(\|v-u\|_{\tilde{L}^{\rho}_{T}(B^{s_{1},
s_{2}}_{p,r})}
\\[2mm]
&\hspace{0,5cm}(1+\|u\|_{L^{\infty}_{T}(L^{\infty})}+\|v\|_{L^{\infty}_{T}(L^{\infty})})+
\|v-u\|_{L^{\infty}_{T}(L^{\infty})}(\|u\|_{\tilde{L}^{\rho}_{T}(B^{s_{1},s_{2}}_{p,r})}+
\|v\|_{\tilde{L}^{\rho}_{T}(B^{s_{1},s_{2}}_{p,r})}).\\
\end{aligned}
$$
\end{proposition}
The proof is a adaptation of a theorem by J.Y. Chemin and H. Bahouri
in \cite{5BC}.\\
We end this section by recalling some estimates in Besov spaces for transport and heat equations. For more details, the reader is referred to \cite{5CT} and \cite{5DH}.
\begin{proposition}
\label{5transport}
Let $(p,r)\in[1,+\infty]^{2}$ and $s\in(-\min(\NN,\frac{N}{p^{'}}),\NN+1)$. Let $u$ be a vector field such that $\n u$ belongs to
$L^{1}(0,T;B^{\NN}_{p,r}\cap L^{\infty})$. Suppose that $q_{0}\in B^{s}_{p,r}$, $F\in L^{1}(0,T, B^{s}_{p,r})$ and that $q\in
L^{\infty}_{T}(B^{s}_{p,r})\cap
C([0,T];{\cal S}^{'})$ solves the following transport equation:
$$
\begin{cases}
\begin{aligned}
&\p_{t}q+u\cdot\n q=F,\\
&q_{t=0}=q_{0}.\\
\end{aligned}
\end{cases}
$$
Let $U(t)=\int^{t}_{0}\|\n u(\tau)\|_{B^{\NN}_{p,r}\cap L^{\infty}}d\tau$. There
exits a constant $C$ depending only on $s$, $p$ and $N$, and such that
for all $t\in [0,T]$, the following inequality holds:
$$\|q\|_{\widetilde{L}^{\infty}_{t}(B^{s}_{p,r})}\leq\exp^{C U(t)}\biggl(\|q_{0}\|_{B^{s}_{p,r}}
+\int^{t}_{0}\exp^{-C U(\tau)}\|F(\tau)\|_{B^{s}_{p,r}}\,d\tau\biggl)$$
If $r<+\infty$ then $q$ belongs to $C([0,T];B^{s}_{p,r})$.
\end{proposition}
\label{5chaleur}
Actually, in \cite{5DH}, the proposition below is proved for non-homogeneous Besov spaces. The adaptation to homogeneous spaces is straightforward. Let us now some estimates for the heat equation:
\begin{proposition}
\label{5chaleur} Let $s\in\R$, $(p,r)\in[1,+\infty]^{2}$ and
$1\leq\rho_{2}\leq\rho_{1}\leq+\infty$. Assume that $u_{0}\in B^{s}_{p,r}$ and $f\in\widetilde{L}^{\rho_{2}}_{T}
(\widetilde{B}^{s-2+2/\rho_{2}}_{p,r})$.
Let u be a solution of:
$$
\begin{cases}
\begin{aligned}
&\p_{t}u-\mu\D u=f\\
&u_{t=0}=u_{0}\,.\\
\end{aligned}
\end{cases}
$$
Then there exists $C>0$ depending only on $N,\mu,\rho_{1}$ and
$\rho_{2}$ such that:
$$\|u\|_{\widetilde{L}^{\rho_{1}}_{T}(\widetilde{B}^{s+2/\rho_{1}}_{p,r})}\leq C\big(
 \|u_{0}\|_{B^{s}_{p,r}}+\mu^{\frac{1}{\rho_{2}}-1}\|f\|_{\widetilde{L}^{\rho_{2}}_{T}
 (\widetilde{B}^{s-2+2/\rho_{2}}_{p})}\big)\,.$$
 If in addition $r$ is finite then $u$ belongs to $C([0,T],B^{s}_{p,r})$.
\end{proposition}
The proof of local well-posedness for initial density bounded away
from zero requires estimates in $B^{s}_{p}$ spaces for the following
linear system:
$$
\begin{cases}
\begin{aligned}
&\p_{t}u-\bar{\mu}{\rm div}(a\n u)-(\bar{\lambda}+\bar{\mu})\n(a{\rm div}u)=G,\\
&u_{t=0}=u_{0},\\
\end{aligned}
\label{512}
\end{cases}
$$
where $u(t,x)\in \R^{N}$, $\bar{\nu}=2\bar{\mu}+\bar{\lambda}>0$ and
the diffusion coefficient is assumed to satisfy:
\begin{equation}
0<\underline{a}\leq a(t,x)\leq\overline{a}. \label{5hip}
\end{equation}
We can prove that the solution of the previous system satisfy
estimates analogous to those of proposition \ref{5chaleur}, see \cite{5DG}.
\begin{proposition}
\label{5chaleurgene}
Let $1<p<+\infty$, $1\leq\alpha_{1}\leq r\leq +\infty$ and s be such
that $\max (1,\NN)\leq s\leq\NN+1$. Set
$\alpha_{2}^{'}=\frac{2}{s-\NN}$ and let $\alpha_{2}$ be such that
$\frac{1}{\alpha_{1}}=\frac{1}{\alpha_{2}}
+\frac{1}{\alpha_{2}^{'}}$. Let assumption (\ref{5hip}) be fulfilled
and $u$ be a solution of (\ref{512}). We suppose that the regularity
index $\tau$ satisfies:
$$1-\frac{2}{\alpha_{1}}-\NN<\tau\leq 1-\frac{2}{\alpha_{1}}+s.$$
Then the following estimate holds for all $\alpha\in[r,+\infty]$:
$$\|u\|_{\widetilde{L}^{\alpha}_{T}(B_{p}^{\tau+\frac{2}{\alpha}})}\lesssim\|u_{0}\|_{B_{p}^{\tau}}+
\|G\|_{\widetilde{L}^{r}_{T}(B_{p}^{\tau+\frac{2}{r}-2})}+\|\n
a\|_{\widetilde{L}^{\alpha^{'}_{2}} _{T}(B_{p}^{s-1})}\|\n
u\|_{\widetilde{L}^{\alpha^{'}_{2}}
_{T}(B_{p}^{\tau-1+\frac{2}{\alpha_{2}}})}.$$
\end{proposition}
\section{Proof of theorem \ref{5theo1}}
\label{5section4}
\subsection{Sketch of the Proof}
In this section, we give the sketch of the proof of theorem
\ref{5theo1} on the global existence result with small initial
data.\\
We will suppose that $\rho$ is close to a constant state
$\bar{\rho}$, so that $\rho$ will be strictly superior to a positive
constant, we will use the parabolicity of the momentum equation to
get a gain of derivatives on the velocity $u$. The density $\rho$
has a behavior similar to the solution of a transport equation. Let
us rewrite the $(SW)$ system in a non conservative form in using the
definition \ref{5bar}.
$$
\begin{cases}
\begin{aligned}
&\p_{t}q+u.\n q+{\rm div}u=F\\[2mm]
&\p_{t}u+u.\n u-\frac{\mu(\bar{\rho})}{\bar{\rho}}\D
u-\frac{\mu(\bar{\rho})+\lambda(\bar{\rho})}{\bar{\rho}}\n{\rm div}u
+(\kappa\bar{\rho}+ P^{'}(\bar{\rho}))\,\n q\\
&\hspace{9,5cm}-\kappa\bar{\rho}\phi*\n q=G\\
\end{aligned}
\end{cases}
\leqno{(SW1)}
$$
where we have:
$$
\begin{aligned}
&F=-q\,{\rm div}u,\\[2mm]
&G={\cal A}(\rho,u)+K(\rho)\n q\,,\\
\end{aligned}
$$
where we set:
$$
\begin{aligned}
&{\cal A}(\rho,u)=\big[\frac{{\rm
div}\big(\mu(\rho)D(u)\big)}{\rho}-\frac{\mu(\bar{\rho})}{\bar{\rho}}\D
u] +\big[\frac{\n\big((\mu(\rho)+\lambda(\rho)){\rm
div}u\big)}{\rho}-\frac{\mu(\bar{\rho})+
\lambda(\bar{\rho})}{\bar{\rho}}\n{\rm div}u\big],\\
&K(\rho)=\frac{\bar{\rho}\,P^{'}(\rho)}{\rho}-P^{'}(\bar{\rho}).\\
\end{aligned}
$$
For $s\in\R$, we denote
$\Lambda^{s}h={\cal{F}}^{-1}(|\xi|^{s}\widehat{h})$.
We set now:
$$d=\Lambda^{-1}{\rm div}u\;\;\mbox{and}\;\;\Omega=\Lambda^{-1}{\rm curl} u$$
where $d$ represents the compressible part of the velocity and
$\Omega$ the incompressible part. We rewrite now the system $(SW1)$
 in using these previous notations on a linear form:
$$
\begin{cases}
\begin{aligned}
&\p_{t}q+\Lambda d=F_{1},\\
&\p_{t}d-\bar{\nu}\D d-\bar{\delta}\Lambda q+\bar{\kappa}\La(\phi*q)=G_{1}\\
&\p_{t}\Omega-\bar{\mu}\D \Omega=H_{1}\\
&u=-\Lambda^{-1}\n d-\La{\rm div}\Omega\\
\end{aligned}
\end{cases}
\leqno{(SW2)}
$$
where we have :
$$\bar{\mu}=\frac{\mu(\bar{\rho})}{\bar{\rho}},\;\;\;\bar{\lambda}=\frac{\lambda(\bar{\rho})}{\bar{\rho}},\;\;\;
\bar{\nu}=2\bar{\mu}+\bar{\lambda},\;\;\;\bar{\delta}=\kappa\bar{\rho}+
P^{'}(\bar{\rho}),\;
\;\;\mbox{and}\;\;\;\bar{\kappa}=\kappa\bar{\rho}.
$$
We have in our case:
$$
\begin{aligned}
&F_{1}=-q\,{\rm div}u-u\cdot\n q,\\
&G_{1}=-\La^{-1}{\rm div}(G),\\
&H_{1}=-\La^{-1}{\rm curl}(G).\\
\end{aligned}
$$
The first idea will be to study the linear system associated to
$(SW2)$. We concentrate on the first two equations because the third
equation is just a heat equation with a non linear term. The system
we want to study reads:
$$
\begin{cases}
\begin{aligned}
&\p_{t}q+\Lambda d=F^{'},\\
&\p_{t}d-\bar{\nu}\D d-\bar{\delta}\Lambda q+\bar{\kappa}\La(\phi*q)=G^{'}.\\
\end{aligned}
\end{cases}
$$
This system has been studied by D. Hoff and K. Zumbrum in \cite{5HZ}
in the case $\bar{\kappa}=0$. There, they investigate the decay
estimates, and exhibit the parabolic smoothing effect on $d$ and on
the low frequencies of
$q$, and a damping effect on the high frequencies of $q$.\\
The problem is that if we focus on this linear system, it appears
impossible to control the term of convection $u\cdot\n q$ which is
one derivative less regular than $q$. Hence we shall include the
convection term in the linear system. We thus have to study:
$$
\begin{cases}
\begin{aligned}
&\p_{t}q+v\cdot\n q+\Lambda d=F,\\
&\p_{t}d+v\cdot\n d-\bar{\nu}\D d-\bar{\delta}\Lambda q+\bar{\kappa}\La(\phi*q)=G,\\
\end{aligned}
\end{cases}
\leqno{(SW2)^{'}}
$$
where $v$ is a function  and we will precise its regularity in the
next proposition. System $(SW2)^{'}$ has been
studied in the case where $\phi=0$ by R. Danchin in \cite{5DG}, we
then adapt the proof in taking into consideration the term coming
from the
capillarity.\\
In the sequel we will assume $\bar{\nu}>0$ and
$\bar{\delta}-\bar{\kappa}\|\widehat{\phi}\|_{L^{\infty}}\geq c>0$.
We get then the following proposition.
\begin{proposition}
\label{5linear1} Let $(q,d)$ a solution of the system $(SW2)^{'}$ on
$[0,T[$ , $1-\N<s\leq 1+\N$ and
$V(t)=\int^{t}_{0}\|v(\tau)\|_{B^{\N+1}}d\tau$. We have then the
following estimate:
$$
\begin{aligned}
&\|(q,d)\|_{\widetilde{B}^{s-1,s}\times B^{s-1}}+\int^{t}_{0}\|(q,d)(\tau)\|_{\widetilde{B}^{s+1,s}\times B^{s+1}}d\tau\\
&\hspace{2cm}\leq Ce^{CV(t)}\big(\|
(q_{0},d_{0})\|_{\widetilde{B}^{s-1,s}\times B^{s-1}}+\int^{t}_{0}
e^{-CV(\tau)}\|
(F,G)(\tau)\|_{\widetilde{B}^{s-1,s}\times B^{s-1}}d\tau\big),\\
\end{aligned}
$$
where $C$ depends only on $\bar{\nu}$, $\bar{\delta}$,
$\bar{\kappa}$, $\phi$, $N$ and $s$.
\end{proposition}
{\bf Proof of Proposition \ref{5linear1}:}\\
\\
Let $(q,u)$ be a solution of $(SW2)^{'}$ and we set:
\begin{equation}
\widetilde{q}=e^{-KV(t)}q,\;\widetilde{u}=e^{-KV(t)}u,\;\widetilde{F}=e^{-KV(t)}F\;\;
\mbox{and}\;\;\widetilde{G}=e^{-KV(t)}G. \label{5def}
\end{equation}
We are going to separate the case of the low and high frequencies,
which have a
different behavior concerning the control of the derivative index for the Besov spaces.\\
In this goal we will consider the two different expressions in low
and high frequencies where $l_{0}\in\mathbb{Z}$, $A$, $B$ and
$K_{1}$ will be fixed later in the proof:
\begin{equation}
\begin{aligned}
&f_{l}^{2}=\bar{\de}\|\widetilde{q}_{l}\|_{L^{2}}^{2}-\ka(\widetilde{q}_{l},\phi*\widetilde{q}_{l})+\|\widetilde{d}_{l}\|_{L^{2}}^{2}-2K_{1}(\Lambda\widetilde{q}_{l},\widetilde{d}_{l})
\;\;\;\;\mbox{for}\;\;l\leq l_{0},\\[2mm]
&f_{l}^{2}=\|\Lambda\widetilde{q}_{l}\|_{L^{2}}^{2}+A\|\widetilde{d}_{l}\|_{L^{2}}^{2}-
\frac{2}{\bar{\nu}}(\La\widetilde{q}_{l},\widetilde{d}_{l})
\;\;\;\;\mbox{for}\;\;l>l_{0}.\\
\end{aligned}
\label{5fl}
\end{equation}
In the first two steps, we show that $K_{1}$ and $A$ may be chosen
such that:
\begin{equation}
2^{l(s-1)}f_{l}^{2}\approx
2^{ls}\max(1,2^{-l})\|\widetilde{q}_{l}\|_{L^{2}}^{2}+2^{l(s-1)}\|
\widetilde{d}_{l}\|_{L^{2}}^{2}, \label{5fl1}
\end{equation}
and we will show the following inequality:
\begin{equation}
\begin{aligned}
\frac{1}{2}\frac{d}{dt}f_{l}^{2}+\alpha\min(2^{2l},1)f_{l}^{2}\leq
C2^{-l(s-1)}\alpha_{l}f_{l}&\big(\,
\|(\widetilde{F},\widetilde{G})\|_{\widetilde{B}^{s-1,s}\times
B^{s-1}}\\[2mm]
&+V^{'}
\|(\widetilde{q},\widetilde{d})\|_{\widetilde{B}^{s-1,s}\times B^{s-1}}\big)
-KV^{'}f_{l}^{2}.\\
\end{aligned}
\label{5inegal}
\end{equation}
where $\sum_{l\in\mathbb{Z}}\alpha_{l}\leq1$ and $\alpha$ is a positive constant.\\
This inequality enables us to get a decay for $q$ and $d$ which will
be used to show a smoothing parabolic effect on $d$.
\subsubsection*{Case of low frequencies}
Applying operator $\D_{l}$ to the system $(SW2)^{'}$, we obtain then
in setting:
$$\widetilde{q}_{l}=\Delta_{l}\widetilde{q},\;\widetilde{d}_{l}=\Delta_{l}\widetilde{d}.$$
the following system:
\begin{equation}
\begin{cases}
\begin{aligned}
&\frac{d}{dt}\widetilde{q}_{l}+\D_{l}(v\cdot\n\widetilde{q})+\Lambda
\widetilde{d}_{l}=\widetilde{F}_{l}-KV^{'}(t)
\widetilde{q}_{l},\\[2mm]
&\frac{d}{dt}d_{l}+\D_{l}(v\cdot\n\widetilde{d}_{l})-\bar{\nu}\Delta
\widetilde{d}_{l}-\bar{\delta}\Lambda\widetilde{q}_{l}+\ka\Lambda(\phi*\widetilde{q}_{l})=\widetilde{G}_{l}-KV^{'}(t)\widetilde{d}_{l}.\\
\end{aligned}
\end{cases}
\label{5LPH}
\end{equation}
We set:
\begin{equation}
f_{l}^{2}=\bar{\de}\|\widetilde{q}_{l}\|_{L^{2}}^{2} +\|\widetilde{d}_{l}
\|_{L^{2}}^{2}-2K_{1}(\Lambda \widetilde{q}_{l},\widetilde{d}_{l})
\label{5P0}
\end{equation}
for some $K_{1}\geq0$  to be fixed hereafter and $(\cdot,\cdot)$
noting the $L^{2}$ inner product.
\\
To begin with, we consider the case  where $F=G=0$, $v=0$ and $K=0$.
Taking the $L^{2}$ scalar product of the first equation of
(\ref{5LPH}) with $\tilde{q}_{l}$ and of the second equation with
$\tilde{d}_{l}$, we get the following two identities:
\begin{equation}
\begin{cases}
\begin{aligned}
&\frac{1}{2}\frac{d}{dt}\|q_{l}\|_{L^{2}}^{2}+(\Lambda d_{l},q_{l})=0,\\[2mm]
&\frac{1}{2}\frac{d}{dt}\|d_{l}\|_{L^{2}}^{2}+\bar{\nu}\|\Lambda
d_{l}\|_{L^{2}}^{2} -\del(\Lambda q_{l},d_{l})+\ka(\Lambda
(\phi*q_{l}),d_{l})
=0.\\
\end{aligned}
\end{cases}
\label{5P1}
\end{equation}
In the same way we have:
\begin{equation}
\begin{aligned}
&\frac{1}{2}\frac{d}{dt}(q_{l},q_{l}*\phi)+(\Lambda d_{l},\phi*q_{l})=0,\\[2mm]
\end{aligned}
\label{5HAS2}
\end{equation}
because we have by the theorem of Plancherel:
$$(\frac{d}{dt}q_{l},q_{l}*\phi)=(\frac{d}{dt}\widehat{q_{l}},\widehat{q_{l}}\widehat{\phi})=
\frac{1}{2}\frac{d}{dt}(\widehat{q_{l}},\widehat{q_{l}}\widehat{\phi})=\frac{1}{2}\frac{d}{dt}(q_{l},q_{l}*\phi).$$
We want now get an equality involving $\bar{\nu}(\Lambda
d_{l},q_{l})$. To achieve it, we apply $\bar{\nu}\Lambda$ to the
first equation of (\ref{5LPH}) and take the $L^{2}$-scalar product
with $d_{l}$, then take the scalar product of the second equation
with $\Lambda q_{l}$ and sum both equalities, which yields:
\begin{equation}
\begin{aligned}
&\frac{d}{dt}(\Lambda q_{l},d_{l})+\|\Lambda d_{l}\|_{L^{2}}^{2}-
\del\|\Lambda q_{l}\|_{L^{2}}^{2}+\ka\|\phi*\La
q_{l}\|_{L^{2}}^{2}+\bar{\nu}(\Lambda^{2}d_{l}, \Lambda q_{l})
=0.\\
\end{aligned}
\label{5P2}
\end{equation}
By linear combination of (\ref{5P1}) and (\ref{5P2}), we get:
\begin{equation}
\begin{aligned}
\frac{1}{2}\frac{d}{dt}f_{l}^{2}+(\bar{\nu}-K_{1})\|\Lambda
d_{l}\|_{L^{2}}^{2} +K_{1}(\del\|\Lambda
q_{l}\|_{L^{2}}^{2}-\ka\|\phi*\La
q_{l}\|_{L^{2}}^{2})-\bar{\nu}K_{1}(\Lambda^{2}d_{l},
\Lambda q_{l})=0.\\
\end{aligned}
\label{5P3}
\end{equation}
And as we have assumed that:
$\de-\ka\|\widehat{\phi}\|_{L^{\infty}}\geq c>0$ we get:
\begin{equation}
\begin{aligned}
\frac{1}{2}\frac{d}{dt}f_{l}^{2}+(\bar{\nu}-K_{1})\|\Lambda
d_{l}\|_{L^{2}}^{2}
+K_{1}c\|q_{l}\|_{L^{2}}^{2}-\bar{\nu}K_{1}(\Lambda^{2}d_{l},
\Lambda q_{l})\leq0.\\
\end{aligned}
\label{5Has3}
\end{equation}
Using spectral localization for $d_{l}$ and convex inequalities, we
find for every $a>0$:
$$
\begin{aligned}
&|(\Lambda^{2}d_{l},\Lambda
q_{l})|\leq\frac{a2^{2l_{0}}}{2}\|\Lambda d_{l}\|_{L^{2}}^{2}
+\frac{1}{2a}\|\Lambda q_{l}\|_{L^{2}}^{2}.\\
\end{aligned}
$$
In using the previous inequality and (\ref{5P3}), we get:
\begin{equation}
\frac{1}{2}\frac{d}{dt}f_{l}^{2}+(\bar{\nu}-K_{1}-\frac{a2^{2l_{0}}}{2})\|\Lambda
d_{l}\|^{2}_{L^{2}} +(K_{1}c-\frac{1}{2a})\|\Lambda
q_{l}\|^{2}_{L^{2}}\leq 0. \label{5P4}
\end{equation}
From (\ref{5P0}) and (\ref{5P4}) we get in choosing $a=\bar{\nu}$
and $K_{1}<\min(\frac{1}{2^{2l_{0}}},
\frac{\bar{\nu}}{2+2^{2l_{0}}\bar{\nu}^{2}})$, then:
\begin{equation}
\frac{1}{2}\frac{d}{dt}f_{l}^{2}+\alpha 2^{2l}f_{l}^{2}\leq 0,
\label{5P5}
\end{equation}
for a constant $\alpha$ depending only on $\bar{\nu}$ and $K_{1}$.\\
In the general case where $F$, $G$, $K$ and $v$ are not zero, we
have:
$$
\begin{aligned}
&\frac{1}{2}\frac{d}{dt}f_{l}^{2}+(\alpha
2^{2l}+KV^{'})f_{l}^{2}\leq
(\widetilde{F}_{l},\widetilde{q}_{l})+(\widetilde{G}_{l},\widetilde{d}_{l})
-K(\Lambda\widetilde{F}_{l},\widetilde{d}_{l})-K(\Lambda\widetilde{G}_{l},\widetilde{q}_{l})
-(\D_{l}(v\cdot\n\widetilde{q}),\widetilde{q}_{l})\\[2mm]
&\hspace{4,5cm}-(\D_{l}(v\cdot\n\widetilde{d}),\widetilde{d}_{l})+K
\big((\Lambda\D_{l}(v\cdot\n\widetilde{q}),\widetilde{d}_{l}\big)
+\big((\Lambda\D_{l}(v\cdot\n\widetilde{d}),\widetilde{q}_{l}\big).\\
\end{aligned}
$$
Now we can use a lemma of harmonic analysis in \cite{5DG} to
estimate the last terms, and get the existence of a sequence
$(\alpha_{l})_{l\in\mathbb{Z}}$ such that
$\sum_{l\in\mathbb{Z}}\alpha_{l}\leq1$ and:
\begin{equation}
\begin{aligned}
&\frac{1}{2}\frac{d}{dt}f_{l}^{2}+(\alpha
2^{2l}+KV^{'})f_{l}^{2}\lesssim
\alpha_{l}f_{l}2^{-l(s-1)}\big(\|(\widetilde{F},\widetilde{G})\|_{\widetilde{B}^{s-1,s}\times
B^{s-1}}
+V^{'}\|(\widetilde{q},\widetilde{d})\|_{\widetilde{B}^{s-1,s}\times B^{s-1}}\big).\\
\end{aligned}
\label{5P6}
\end{equation}
\subsubsection*{Case of high frequencies}
We consider now the case where $l\geq l_{0}+1$ and we recall that:
$$f_{l}^{2}=\|\Lambda\widetilde{q}_{l}\|_{L^{2}}^{2}+A\|\widetilde{d}_{l}\|_{L^{2}}^{2}-
\frac{2}{\bar{\nu}}(\widetilde{q}_{l},\widetilde{d}_{l}).$$ For the
sake of simplicity, we suppose here that $F=G=0$, $v=0$ and $K=0$.
We now want a control $\|\Lambda q_{l}\|_{L^{2}}^{2}$ on e apply the
operator $\Lambda$ to the first equation of (\ref{5LPH}), multiply
by $\Lambda q_{l}$ and integrate over $\R^{N}$, so we obtain:
\begin{equation}
\frac{1}{2}\frac{d}{dt}\|\Lambda
q_{l}\|_{L^{2}}^{2}+(\Lambda^{2}d_{l},\Lambda q_{l})=0.
 \label{5H7}
\end{equation}
Moreover we have:
\begin{equation}
\begin{aligned}
&\frac{1}{2}\frac{d}{dt}\|d_{l}\|_{L^{2}}^{2}+\bar{\nu}\|\Lambda
d_{l}\|_{L^{2}}^{2} -\del(\Lambda q_{l},d_{l})+\ka(\Lambda
(\phi*q_{l}),d_{l})
=0.\\[2mm]
&\frac{d}{dt}(\Lambda q_{l},d_{l})+\|\Lambda d_{l}\|_{L^{2}}^{2}-
\del\|\Lambda q_{l}\|_{L^{2}}^{2}+\ka\|\phi*\La
q_{l}\|_{L^{2}}^{2}+\bar{\nu}(\Lambda^{2}d_{l}, \Lambda q_{l})
=0.\\
\end{aligned}
\label{5H8}
\end{equation}
By linear combination of (\ref{5H7})-(\ref{5H8}) we have:
\begin{equation}
\begin{aligned}
&\frac{1}{2}\frac{d}{dt}f_{l}^{2}+\frac{1}{\bar{\nu}}\|\Lambda
q_{l}\|_{L^{2}}^{2}+\big(A\bar{\nu}-\frac{1}{\bar{\nu}}\big)\|\Lambda
d_{l}\|_{L^{2}}^{2}-A\del(\Lambda
q_{l},d_{l})+A\ka(\Lambda (\phi*q_{l}),d_{l})=0.\\
\end{aligned}
\label{5P10}
\end{equation}
Moreover we have:
$$|-A\del(\Lambda
q_{l},d_{l})+A\ka(\Lambda (\phi*q_{l}),d_{l})|\leq
A(\del+\ka\|\widehat{\phi}\|_{L^{\infty}})|(\Lambda q_{l},d_{l})|$$
We have now in using Young inequalities for all $a>0$:
$$
\begin{aligned}
&|(d_{l},\Lambda q_{l})|\leq\frac{a}{2}\|\Lambda q_{l}\|_{L^{2}}^{2}+\frac{1}{2a}\|d_{l}\|_{L^{2}}^{2},\\
\end{aligned}
$$
So we get:
\begin{equation}
\begin{aligned}
&\frac{1}{2}\frac{d}{dt}f_{l}^{2}+2^{2l_{0}}\big(A\bar{\nu}-\frac{1}{\bar{\nu}}-\frac{1}{2a}\big)
\|d_{l}\|_{L^{2}}^{2}+(\frac{1}{\bar{\nu}}-\frac{a}{2})\|\Lambda
q_{l}\|_{L^{2}}^{2}.
\leq 0\\
\end{aligned}
\label{5P12}
\end{equation}
So in choosing:
$$a=\frac{1}{\bar{\nu}A}\;\;\;\mbox{and}\;\;\;A>\max(\frac{2}{\bar{\nu}},1)$$
there exists a constant $\alpha$ such that for $l\geq l_{0}+1$ we
have:
\begin{equation}
\frac{1}{2}\frac{d}{dt}f_{l}^{2}+\alpha f_{l}^{2}\leq0. \label{5P13}
\end{equation}
In the general case where $F$, $G$, $H$, $K$ and $v$ are not
necessarily zero, we use a lemma of harmonic analysis in \cite{5DG}
to control the convection terms. We finally get:
\begin{equation}
\begin{aligned}
&\frac{1}{2}\frac{d}{dt}f_{l}^{2}+(\alpha+KV^{'})f_{l}^{2}\,\lesssim
\alpha_{l}f_{l}2^{-l(s-1)}\big(\|(\widetilde{F},\widetilde{G})\|_{\widetilde{B}^{s-1,s}\times B^{s-1}}\\[2mm]
&\hspace{8,5cm}+V^{'}\|(\widetilde{q},\widetilde{d})\|_{\widetilde{B}^{s-1,s}\times B^{s-1}}\big).\\
\end{aligned}
\end{equation}
This finish the proof of (\ref{5fl}) and (\ref{5inegal}).
\subsubsection*{The damping effect}
We are now going to show that inequality (\ref{5inegal}) entails a
decay for $q$ and $d$.
In fact we get a parabolic decay for $d$, while $q$ has a behavior similar to a transport equation.\\
Using $h_{l}^{2}=f_{l}^{2}+\delta^{2}$, integrating over $[0,t]$ and
then having $\delta$ tend to 0, we infer:
\begin{equation}
\begin{aligned}
f_{l}(t)+\alpha&\min(2^{2l},1)\int^{t}_{0}f_{l}(\tau)d\tau\\[2mm]
&\leq
f_{l}(0)+C2^{-l(s-1)}\int^{t}_{0}\alpha_{l}(\tau)\|(\widetilde{F}(\tau),\widetilde{G}(\tau))
\|_{\widetilde{B}^{s-1,s}\times B^{s}}d\tau\\[2mm]
&\hspace{2cm}+\int^{t}_{0}V^{'}(\tau)\big(C2^{-l(s-1)}\alpha_{l}(\tau)\|
(\widetilde{q},\widetilde{d})\|_{\tilde{B}^{s-1,s}\times B^{s}}-Kf_{l}(\tau)\big)d\tau.\\
\end{aligned}
\label{5P14}
\end{equation}
Thanks to (\ref{5fl1}), we have in taking $K$ large enough :
$$\sum_{l\in\mathbb{Z}}\big(C2^{-l(s-1)}\alpha_{l}(\tau)\|
(\widetilde{q},\widetilde{d})\|_{\tilde{B}^{s-1,s}\times B^{s}}
-Kf_{l}(\tau)\big)\leq0,$$
In multiplying (\ref{5P14}) by $2^{l(s-1)}$ and in using the last
inequality, we conclude after summation on $\mathbb{Z}$, that:
\begin{equation}
\begin{aligned}
&\|\widetilde{q}(t)\|_{\tilde{B}^{s-1,s}}+\|\widetilde{d}\|_{\tilde{B}^{s-1}}
+\alpha\int^{t}_{0}\|\widetilde{q}(\tau)\|_{\widetilde{B}^{s-1,s}}d\tau+\sum_{l\in\mathbb{Z}}\int^{t}_{0}\alpha2^{l(s-1)}
\min(2^{2l},1)\|\widetilde{d}_{l}(\tau)\|_{L^{2}}d\tau\\[2mm]
&\hspace{6cm}\lesssim
\|(\widetilde{q}_{0},\widetilde{d}_{0})\|_{\widetilde{B}^{s-1,s}\times
B^{s-1}} +\int^{t}_{0}
\|(\widetilde{F},\widetilde{G})\|_{\widetilde{B}^{s-1,s}\times B^{s-1}}d\tau.\\
\end{aligned}
\label{5P15}
\end{equation}
\subsubsection*{The smoothing effect}
Once stated the damping effect for $q$, it is easy to get the
smoothing effect on $d$ by considering the last two equations where
the term $\Lambda q$ is considered as a source term
.\\
Thanks to (\ref{5P15}), it suffices to prove it for high frequencies
only. We therefore suppose in this subsection that
$l\geq l_{0}$ for a $l_{0}$ big enough.\\
We set $g_{l}=\|\widetilde{d}_{l}\|_{L^{2}}$
and in using the previous inequalities, we have:
$$\frac{1}{2}\frac{d}{dt}\|\widetilde{d}_{l}\|_{L^{2}}^{2}+\bar{\nu}\|\Lambda \widetilde{d}_{l}\|_{L^{2}}^{2}
-\del(\Lambda \widetilde{q}_{l},\widetilde{d}_{l})+\ka(\Lambda
(\phi*\widetilde{q}_{l}),\widetilde{d}_{l}) =\widetilde{G}_{l}\cdot
\widetilde{d}_{l}-KV^{'}(t)\|\widetilde{d}_{l}\|_{L^{2}}^{2}.$$ We
get finally with $\alpha>0$:
$$
\begin{aligned}
\frac{1}{2}\frac{d}{dt}g_{l}^{2}+\alpha2^{2l}g_{l}^{2}\leq\,
g_{l}\big(\|\Lambda\widetilde{q}_{l}\|_{L^{2}}&+\|\widetilde{G}_{l}\|_{L^{2}}\big)
+g_{l}V^{'}(t)(C\alpha_{l}2^{-l(s-1)}\|\widetilde{d}\|_{B^{s-1}}-Kg_{l}\big).\\
\end{aligned}
$$
\\
We therefore get in using standard computations:
$$
\begin{aligned}
&\sum_{l\geq l_{0}}2^{l(s-1)}\|\widetilde{d}_{l}(t)\|_{L^{2}}+
\alpha\int^{t}_{0}\sum_{l\geq
l_{0}}2^{l(s+1)}\|\widetilde{d}_{l}(\tau)\|_{L^{2}}
d\tau\leq\|d_{0}\|_{B^{s-1}}
+\int^{t}_{0}\|\widetilde{G}(\tau)\|_{B^{s-1}}
d\tau\\[2mm]
&\hspace{5,5cm}+\int^{t}_{0}\sum_{l\geq
l_{0}}2^{ls}\|\widetilde{q}_{l}(\tau)\|_{L^{2}}+CV(t)\sup_{\tau\in[0,t]}(\|\widetilde{d}(\tau)\|_{B^{s-1}}).\\
\end{aligned}
$$
Using the above inequality and (\ref{5P15}), we have:
\begin{equation}
\begin{aligned}
&\int^{t}_{0}\sum_{l\geq
l_{0}}2^{l(s+1)}\|\widetilde{d}_{l}(\tau)\|_{L^{2}} d\tau\lesssim
(1+V(t))\big(\|q_{0}\|_{\widetilde{B}^{s-1,s}}+\|d_{0}\|_{B^{s-1}}
\big)\\[2mm]
&\hspace{7cm}+\int^{t}_{0}(\|\widetilde{F}(\tau)\|_{\widetilde{B}^{s-1,s}}+\|\widetilde{G}(\tau)\|_{B^{s-1}}
)d\tau.\\[2mm]
\end{aligned}
\label{5P16}
\end{equation}
Combining that last inequality (\ref{5P16}) with (\ref{5P15}), we
achieve the proof of
proposition \ref{5linear1}.\\
\hfill {$\Box$}
\subsection{Proof of theorem \ref{5theo1}}
This section is devoted to the proof of the theorem \ref{5theo1}.
The principle of the proof is a very classical one. We want to
construct a sequence $(q^{n},u^{n})_{n\in\mathbb{N}}$ of approximate
solutions of the system $(SW)$, and we will use the proposition
\ref{5linear1} to get some uniform bounds on
$(q^{n},u^{n})_{n\in\mathbb{N}}$. We will conclude by stating some
properties of compactness, which will guarantee that up to an
extraction, $(q^{n},u^{n})_{n\in\mathbb{N}}$ converges to a solution
$(q,u)$ of the system $(SW)$.
\subsubsection*{First step: Building the sequence $(q^{n},u^{n})_{n\in\mathbb{N}}$}
We start with the construction of the sequence
$(q^{n},u^{n})_{n\in\mathbb{N}}$, in this goal we use the Friedrichs
operators $(J_{n})_{n\in\mathbb{N}}$ defined by:
$$J_{n}g={\cal F}^{-1}(1_{B(\frac{1}{n},n)}\widehat{g}),$$
where ${\cal F}^{-1}$ is the inverse Fourier transform.
Let us consider the approximate system:
\begin{equation}
\begin{cases}
\begin{aligned}
&\p_{t}q^{n}+J_{n}(J_{n}u^{n}\cdot\n J_{n}q^{n})+\Lambda J_{n}d^{n}=F^{n}\\
&\p_{t}d^{n}+J_{n}(J_{n}u^{n}\cdot\n J_{n}d^{n})-\bar{\nu}\D
J_{n}d^{n}-\del\Lambda J_{n}q^{n}-\ka
\phi*\Lambda J_{n}q^{n}=G^{n}\\
&\p_{t}\Omega^{n}-\bar{\nu}\D J_{n}\Omega^{n}=H^{n}\\
&u^{n}=-\Lambda^{-1}\n d^{n}-\Lambda^{-1}{\rm div}\Omega^{n}\\
&(q^{n},d^{n},\Omega^{n})_{/t=0}=(J_{n}q_{0},J_{n}d_{0},J_{n}\Omega_{0})\\
\end{aligned}
\end{cases}
\label{5grande equation}
\end{equation}
with:
$$
\begin{aligned}
&F^{n}=-J_{n}\big((J_{n}q^{n}){\rm div}J_{n}u^{n}\big),\\[2mm]
&G^{n}=J_{n}\Lambda^{-1}{\rm div}\big[{\cal
A}(\va\big(\bar{\rho}(1+J_{n}q^{n})\big),J_{n}u^{n})
+K(\va\big(\bar{\rho}(1+J_{n}q^{n})\big)\n q^{n}\big],
\\[2mm]
&H^{n}=J_{n}\Lambda^{-1}{\rm curl}\big[{\cal
A}(\va\big(\bar{\rho}(1+J_{n}q^{n})\big),J_{n}u^{n})
+K(\va\big(\bar{\rho}(1+J_{n}q^{n})\big)\n q^{n}\big].\\
\end{aligned}
$$
where $\va$ is a smooth function verifying $\va(s)=s$ for $\frac{1}{n}\leq s\leq n$ and $\va\geq\frac{1}{4}$.\\
We want to show that (\ref{5grande equation}) is only an ordinary
differential equation in $L^{2}\times L^{2}\times L^{2}$. We can observe easily that all the source term in
(\ref{5grande equation}) turn out to be continuous in $L^{2}\times
L^{2}\times L^{2}$. As a example, we consider
the term $J_{n}{\cal
A}(\va\big(\bar{\rho}(1+J_{n}q^{n})\big),J_{n}u^{n})$. We have then
by Plancherel theorem:
$$
\begin{aligned}
\|J_{n}\big(\frac{{\rm
div}\big(\mu(\va\big(\bar{\rho}(1+J_{n}q^{n})\big)DJ_{n}u^{n}\big)}{\va\big(\bar{\rho}(1+
J_{n}q^{n})\big)}\big)\|_{L^{2}}
&\leq n\|\mu(\va\big(\bar{\rho}(1+J_{n}q^{n})\big)DJ_{n}u^{n}\|_{L^{2}}\\
&\hspace{3cm}\times\|\frac{1}{\va\big(\bar{\rho}(1+
J_{n}q^{n})\big)}\|_{L^{\infty}},\\
&\leq 4M_{n} n^{2}\|u^{n}\|_{L^{2}}.\\
\end{aligned}
$$
where $M_{n}=\|\mu(\va\big(\bar{\rho}(1+J_{n}q^{n})\|_{L^{\infty}}$.\\
According to the Cauchy-Lipschitz theorem, a unique maximal solution
exists in $C([0,T_{n});L^{2})$ with $T_{n}>0$. Moreover, since
$J_{n}=J_{n}^{2}$ we show that
$(J_{n}q^{n},J_{n}d^{n},J_{n}\Omega^{n})$ is also a solution and
then by uniqueness we get that
$(J_{n}q^{n},J_{n}u^{n})=(q^{n},u^{n})$. This implies that
$(q^{n},d^{n},\Omega^{n})$ is solution of the following system:
\begin{equation}
\begin{cases}
\begin{aligned}
&\p_{t}q^{n}+J_{n}(u^{n}.\n q^{n})+\Lambda d^{n}=F_{1}^{n}\\
&\p_{t}d^{n}+J_{n}(u^{n}.\n d^{n})-\bar{\nu}\D d^{n}-\del\Lambda
q^{n}-\ka
\phi*\Lambda q^{n}=G_{1}^{n}\\
&\p_{t}\Omega^{n}-\bar{\nu}\D \Omega^{n}=H_{1}^{n}\\
&u^{n}=-\Lambda^{-1}\n d^{n}-\Lambda^{-1}{\rm div}\Omega^{n}\\
&(q^{n},d^{n},\Omega^{n})_{/t=0}=(J_{n}q_{0},J_{n}d_{0},J_{n}\Omega_{0})\\
\end{aligned}
\end{cases}
\label{5grande equation1}
\end{equation}
and:
$$
\begin{aligned}
&F_{1}^{n}=-J_{n}\big(q^{n}{\rm div}u^{n}\big),\\[2mm]
&G_{1}^{n}=J_{n}\Lambda^{-1}{\rm div}\big[{\cal
A}(\va\big(\bar{\rho}(1+q^{n})\big),u^{n})
+K(\va\big(\bar{\rho}(1+q^{n})\big)\big],
\\[2mm]
&H_{1}^{n}=J_{n}\Lambda^{-1}{\rm curl}\big[{\cal
A}(\va\big(\bar{\rho}(1+q^{n})\big),u^{n})
+K(\va\big(\bar{\rho}(1+q^{n})\big)\big].\\
\end{aligned}
$$
And the system (\ref{5grande equation1}) is again an ordinary differential
equation in $L^{2}_{n}$ with:
$$L^{2}_{n}=\{g\in L^{2}(\R^{N})/ {\rm supp}\widehat{g}\subset B(\frac{1}{n},n)\}.$$
Due to the Cauchy-Lipschitz theorem again, a unique maximal solution
exists in $C^{1}([0,T^{'}_{n});L^{2}_{n})$ with $T^{'}_{n}\geq T_{n}>0$.
\subsubsection*{Second step: Uniform estimates}
In this part, we want to get uniform estimates independent of $T$ on $\|(q^{n},u^{n})\|_{E_{T}^{\frac{N}{2}}}$
for all $T<T^{'}_{n}$. This will show that $T^{'}_{n}=+\infty$ by Cauchy-Lipchitz because
the norms $\|\cdot\|_{E^{\N}}$ and $L^{2}$ are equivalent on $L^{2}_{n}$.
$E_{T}^{\frac{N}{2}}$
Let us set:
$$
\begin{aligned}
&E(0)=\|q_{0}\|_{\widetilde{B}^{\N-1,\N}}+\|u_{0}\|_{B^{\N}},\\
&E(q,u,t)=\|q\|_{L^{\infty}_{t}(\widetilde{B}^{\frac{N}{2}-1,\frac{N}{2}})}+\|q\|_{L^{\infty}_{t}(B^{\frac{N}{2}-1})}
+\|q\|_{L^{1}_{t}(\widetilde{B}^{\frac{N}{2}+1,\frac{N}{2}})}+\|q\|_{L^{\infty}_{t}(B^{\frac{N}{2}+1})},\\
\end{aligned}
$$
and:
$$\bar{T}_{n}=\sup\{t\in[0,T^{'}_{n}),E(q^{n},u^{n},t)\leq 3CE(0)\}$$
$C$ corresponds to the constant in the proposition
\ref{5linear1} and as $C>1$ we have $3C>1$ so by continuity
we have $\bar{T}_{n}>0$.\\
We are going to prove that $\bar{T}_{n}=T^{'}_{n}$ for all $n\in\mathbb{N}$ and we will conclude
that $\forall n\in\mathbb{N}$ $T^{'}_{n}=+\infty$.
To achieve it, one can use the proposition \ref{5linear1} to the system (\ref{5grande equation1}) to obtain
uniform bounds, so we get in setting $V_{n}(t)=\|u^{n}\|_{L^{1}_{T}(B^{\N+1})}$:
$$
\begin{aligned}
\|(q^{n},u^{n})\|_{E_{T}^{\N}}\leq
C\,e^{CV_{n}(t)}\big(\,\|q_{0}\|_{\tilde{B}^{\N-1,\N}}
&+\|u_{0}\|_{B^{\N}}+\int_{0}^{T}e^{-CV_{n}(\tau)}(\|F_{1}^{n}(\tau)\|_{\widetilde{B}^{\N-1,\N}}\\
&+\|G_{1}^{n}(\tau)\|_{B^{\N-1}}+\|H_{1}^{n}(\tau)\|_{B^{\N-1}})d\tau
.\big)\\
\end{aligned}
$$
Therefore, it is only a matter of proving appropriate estimates for
$F_{1}^{n}$, $G_{1}^{n}$ and $H_{1}^{n}$ in using properties
of continuity on the paraproduct.\\
\\
We estimate now
$\|F_{1}^{n}\|_{L^{1}_{T}(\widetilde{B}^{\frac{N}{2}-1,\frac{N}{2}})}$
in using proposition \ref{5produit} and \ref{5composition}:
$$
\begin{aligned}
\|F_{1}^{n}\|_{L^{1}_{T}(B^{\frac{N}{2}-1,\frac{N}{2}})}&\leq C\|q^{n}\|_{L_{T}^{\infty}(B^{\frac{N}{2}-1,\frac{N}{2}})}\|{\rm
div}u^{n}\|_{L^{1}_{T}(B^{\frac{N}{2}})},\\
\end{aligned}
$$
We now want to estimate $G_{1}^{n}$:
$$
\begin{aligned}
\|{\cal A}(\va(\bar{\rho}(1+q^{n})),u^{n})\|_{L^{1}_{T}(B^{\frac{N}{2}-1})} &\leq C
\|u^{n}\|_{L^{1}_{T}(B^{\frac{N}{2}+1})}\|q^{n}\|_{L_{T}^{\infty}(B^{\frac{N}{2}})}(1+
\|q^{n}\|_{L_{T}^{\infty}(B^{\frac{N}{2}})}),\\
\end{aligned}
$$
We can verify that $K$ fulfills the hypothesis of the proposition
\ref{5composition}, so we get:
$$
\begin{aligned}
\|K(\va(\bar{\rho}(1+q^{n}))\n q^{n}\|_{L^{1}_{T}(B^{\frac{N}{2}-1})}&\leq
C\|q^{n}\|^{2}_{L^{2}_{T}(B^{\frac{N}{2}})}\|q^{n}\|_{L_{T}^{\infty}(\widetilde{B}^{\frac{N}{2}-1,\frac{N}{2}})},\\
\end{aligned}
$$
Moreover we recall that according to proposition \ref{5produit}:
$$\|q^{n}\|_{L^{2}_{T}(B^{\frac{N}{2}})}^{2}\leq\|q^{n}\|_{L_{T}^{\infty}(\widetilde{B}^{\frac{N}{2}-1,\frac{N}{2}})}
\|q^{n}\|_{L_{T}^{1}(\widetilde{B}^{\frac{N}{2}+1,\frac{N}{2}})}.
$$
We proceed similarly to estimate
$\|H_{1}^{n}\|_{L_{T}^{1}(B^{\frac{N}{2}-1})}$ and finally we have:
$$
\begin{aligned}
\|F_{1}^{n}\|_{L_{T}^{1}(B^{\frac{N}{2}-1})}+\|G_{1}^{n}\|_{L_{T}^{1}(B^{\frac{N}{2}-1})}+
\|H_{1}^{n}\|_{L^{1}(B_{T}^{\frac{N}{2}-1})}\leq
2C(E^{2}(q^{n}&,u^{n},T)\\
&\hspace{0,5mm}+E^{3}(q^{n},u^{n},T)),
\end{aligned}
$$
whence:
$$\|(q^{n},u^{n})\|_{E_{T}^{\frac{N}{2}}}\leq
Ce^{C^{2}3E(0)}E(0)(1+18CE(0)(1+3E(0))),$$
We want now to get:
$$e^{3C^{2}E(0)}(1+18CE(0)(1+3E(0)))\leq 2$$
for this it suffices choose $E(0)$ small enough, let $E(0)<\e$
such that:
$$1+18CE(0)(1+3E(0))\leq\frac{3}{2}\;\;\;\mbox{and}\;\;\;e^{3C^{2}E(0)}\leq\frac{4}{3}.$$
So we get $\bar{T}_{n}=T^{'}_{n}$, indeed we have shown that $\forall T$
such that $T<\bar{T}_{n}$:
$$E(q^{n},u^{n},T)\leq2CE(0).$$
Then we have $\bar{T}_{n}=T^{'}_{n}$, because if $\bar{T}_{n}<T^{'}_{n}$ we have seen
that $E(q^{n},u^{n},\bar{T}_{n})\leq2CE(0)$ and so by continuity for $\bar{T}_{n}+\e$ with $\e$ small enough
we obtain again $E(q^{n},u^{n},\bar{T}_{n}+\e)\leq3CE(0)$ and stands in contradiction with the definition of $\bar{T}_{n}$.\\
So if $\bar{T}_{n}=T^{'}_{n}<+\infty$ we have seen that:
$$E(q^{n},u^{n},T^{'}_{n})\leq3CE(0).$$
As $\|q_{n}\|_{L_{T^{'}_{n}}^{\infty}(\widetilde{B}^{\N})}<+\infty$ and $\|u_{n}\|_{L_{T^{'}_{n}}^{\infty}(\widetilde{B}^{\N-1})}<+\infty$, it implies
that
$\|q_{n}\|_{L_{T^{'}_{n}}^{\infty}(L^{2}_{n})}<+\infty$ and $\|u_{n}\|_{L_{T^{'}_{n}}^{\infty}(L^{2}_{n})}<+\infty$,
so by Cauchy-Lipschitz theorem, one may continue the solution beyond $T^{'}_{n}$ which contradicts the definition of $T^{'}_{n}$.\\
Finally the approximate solution $(q^{n},u^{n})_{n\in\mathbb{N}}$ is
global in time.
\subsubsection*{Second step: existence of a solution }
In this part, we shall show that, up to an extraction, the sequence
$(q^{n},u^{n})_{n\in\mathbb{N}}$ converges in
${\cal{D}}^{'}(\R^{+}\times\R^{N})$ to a solution $(q,u)$ of $(SW)$
which has the desired regularity properties. The proof lies on
compactness arguments. To start with, we show that the time first
derivative of $(q^{n},u^{n})$ is uniformly bounded in appropriate
spaces. This enables us to apply Ascoli's theorem and get the
existence of a limit $(q,u)$ for a subsequence. Now, the uniform
bounds of the previous part provide us with additional regularity
and convergence properties so that we may pass to the limit in
the system.\\
It is convenient to split $(q^{n},u^{n})$ into the solution of a
linear system with initial data $(q_{n},u_{n})$ and forcing term,
and the discrepancy to that solution.\\
More precisely, we denote by $(q^{n}_{L},u^{n}_{L})$ the solution
to:\\
\begin{equation}
\begin{aligned}
&\p_{t}q^{n}_{L}+{\rm div}u^{n}_{L}=0\\
&\p_{t}u^{n}_{L}-{\cal{A}}u^{n}_{L}+\n q^{n}_{L}=0\\
&(q^{n}_{L},v^{n}_{L})_{/t=0}=(J_{n}q_{0},J_{n}u_{0})\\
\end{aligned}
\label{5qLuL}
\end{equation}
where: ${\cal A}=\bar{\mu}\D+(\bar{\lambda}+\bar{\mu})\n{\rm div}$
and we set $(\bar{q}^{n},\bar{u}^{n})=(q^{n}-q^{n}_{L},u^{n}-u^{n}_{L})$.\\
Obviously, the definition of $(q^{n}_{L},v^{n}_{L})_{/t=0}$
entails:
$$(q^{n}_{L})_{/t=0}\rightarrow
q_{0}\;\mbox{in}\;\widetilde{B}^{\frac{N}{2}-1,\frac{N}{2}},\;
(u^{n}_{L})_{/t=0}\rightarrow
u_{0}\;\mbox{in}\;\tilde{B}^{\frac{N}{2}-1}.$$
The proposition
\ref{5chaleur} insures that $(q^{n}_{L},u^{n}_{L})$ converges to
the solution $(q_{L},u_{L})$ of the linear system associated to
(\ref{5qLuL}) in $E^{\frac{N}{2}}$. We now have to prove the
convergence of $(\bar{q}^{n},\bar{u}^{n})$. This is of course a
trifle more difficult and requires compactness results. Let us first
state the following lemma.
\begin{lemme}
$(q^{n},u^{n}))_{n\in\mathbb{N}}$ is uniformly bounded
in
$C^{\frac{1}{2}}(\R^{+};B^{\frac{N}{2}-1})\times(C^{\frac{1}{4}}(\R^{+};B^{\frac{N}{2}-\frac{3}{2}}))^{N}$.
\label{5lemmereference}
\end{lemme}
{\bf Proof:}\\
\\
In all the proof, we will note u.b for uniformly bounded.\\
We first prove that $\frac{\p}{\p t}q^{n}$ is u.b in
$L^{2}(\R^{+},B^{\frac{N}{2}-1})$, which yields the desired result
for
$q^{n}$.\\
Let us observe that $q^{n}$  verifies the following
equation
$$\frac{\p}{\p t}q^{n}={\rm div}u^{n}-J_{n}(u^{n}.\n q^{n})-J_{n}(q^{n}{\rm div}
u^{n}).$$
According to the first part, $(u_{n})_{n\in\mathbb{N}}$ is u.b in
$L^{2}(B^{\frac{N}{2}})$, so we can conclude that $\frac{\p}{\p
t}q^{n}$ is u.b in $L^{2}(B^{\frac{N}{2}-1})$. Indeed we have:
$$
\begin{aligned}
&\|J_{n}(q^{n}{\rm div}u^{n})\|_{L^{2}(B^{\frac{N}{2}-1})}\leq\|u^{n}\|_{L^{2}(B^{\frac{N}{2}})}
\|q^{n}\|_{L^{\infty}(B^{\frac{N}{2}})},\\
&\|J_{n}(u^{n}.\n q^{n})\|_{L^{2}(B^{\frac{N}{2}-1})}\leq\|u^{n}\|_{L^{2}(B^{\frac{N}{2}})}\|q^{n}\|_{L^{\infty}(B^{\frac{N}{2}})}.\\
\end{aligned}
$$
And we recall that we use the fact that
$\widetilde{B}^{\frac{N}{2}-1,\frac{N}{2}}\hookrightarrow
B^{\frac{N}{2}}$.\\
Let us prove
now that $\frac{\p}{\p t}d^{n}$ is u.b in
$L^{\frac{4}{3}}(B^{\frac{N}{2}-\frac{3}{2}})+L^{4}(B^{\frac{N}{2}-\frac{3}{2}})$
and that $\p_{t}\Omega^{n}$ is u.b in
$L^{\frac{4}{3}}(B^{\frac{N}{2}-\frac{3}{2}})$ (which gives the
required result for $u^{n}$ in using the relation
$u^{n}=-\Lambda^{-1}\n d^{n}-\Lambda^{-1}{\rm div}
\Omega^{n}$).\\
Let us recall that:
$$
\begin{aligned}
&\frac{\p}{\p t}d^{n}=J_{n}(u^{n}\cdot\n d^{n})+J_{n}\Lambda^{-1}{\rm
div}\big[{\cal A}(\va(\bar{\rho}(1+q^{n})),u^{n})
+J_{n}(K(\va(\bar{\rho}(1+q^{n})))\n q^{n})\big]\\
&\hspace{9,1cm}+\bar{\nu}\D d^{n}+\del\Lambda q^{n}-\ka\phi*\La q^{n},\\[4mm]
&\frac{\p}{\p t}\Omega^{n}=J_{n}\Lambda^{-1}{\rm curl}\big[{\cal
A}(\va(\bar{\rho}(1+q^{n})),u^{n})
+J_{n}(K(\va(\bar{\rho}(1+q^{n}))\n q^{n}))\big]+\bar{\mu}\D\Omega^{n}.\\
\end{aligned}
$$
Results of step one and an interpolation argument yield uniform
bounds for $u^{n}$ in $L^{\infty}(B^{\frac{N}{2}-1})\cap
L^{\frac{4}{3}}(B^{\frac{N}{2}+\frac{1}{2}})$, we  infer in
proceeding as for $\frac{\p}{\p t}q^{n}$ that:
$$
\begin{aligned}
A_{n}=J_{n}(u^{n}\cdot\n d^{n})+J_{n}\Lambda^{-1}{\rm
div}\big[{\cal A}(\va(\bar{\rho}(1+q^{n})),u^{n})
+J_{n}(K(\va(\bar{\rho}(&1+q^{n})))\n q^{n})\big]+\bar{\nu}\D d^{n}\\
&\;\;\;\;\;\mbox{is u.b in}\;L^{\frac{4}{3}}(B^{\frac{N}{2}-\frac{3}{2}}).\\
\end{aligned}
$$
Using the bounds for $q^{n}$ in $L^{2}(B^{\frac{N}{2}})\cap
L^{\infty}( \widetilde{B}^{\frac{N}{2}-1,\frac{N}{2}})$, we get
$q^{n}$ u.b in $L^{4}(B^{\frac{N}{2}-\frac{1}{2}})$ in using
proposition \ref{5produit}. We thus have $J_{n}(K(\va(\bar{\rho}(1+q^{n}))\n
q^{n}$ u.b in $L^{\frac{4}{3}}(B^{\frac{N}{2}-\frac{3}{2}})$.\\
Using the bounds for $u^{n}$ in $L^{\infty}(B^{\frac{N}{2}-1})\cap
L^{\frac{4}{3}}(B^{\frac{N}{2}+\frac{1}{2}})$ we
finally get $A_{n}$ is u.b in
$L^{\frac{4}{3}}(B^{\frac{N}{2}-\frac{3}{2}})$.
To conclude $\phi*\La q^{n}$ is u.b in $L^{4}(B^{\N-\frac{3}{2}})$, so $\frac{\p}{\p t}d^{n}$
is u.b in $L^{\frac{4}{3}}(B^{\frac{N}{2}-\frac{3}{2}})+L^{4}(B^{\frac{N}{2}-\frac{3}{2}})$.
\\
\\
The case of $\frac{\p}{\p t}\Omega^{n}$ goes along the same
lines. As the terms corresponding to
$\Lambda q^{n}$ and $\phi*\Lambda\bar{q}^{n}$ do not appear, we
simply get $\p_{t}\Omega^{n}$ u.b in
$L^{\frac{4}{3}}(B^{\frac{N}{2}-\frac{3}{2}})$.\\
\null \hfill {$\Box$}
\\
We can now turn to the proof of the existence of a solution and
using Ascoli theorem to get strong convergence. We proceed similarly
to the theorem of Aubin-Lions.
\begin{theorem}
Let $X$ a compact metric space and $Y$ a complete metric space. Let
$A$ be an equicontinuous part of $C(X,Y)$. Then we have the two
equivalent proposition:
\begin{enumerate}
\item  $A$ is relatively compact in  $C(X,Y)$
\item $A(x)=\{ f(x);\;\;f\in A\}$ is relatively compact in $Y$
\end{enumerate}
\end{theorem}
We need to localize because we have some result of compactness for
the local Sobolev space. Let $(\chi_{p})_{p\in\mathbb{N}}$ be a
sequence of $C^{\infty}_{0}(\R^{N})$ cut-off functions supported in
the ball $B(0,p+1)$ of $\R^{N}$ and equal to 1 in a neighborhood of
$B(0,p)$.\\
For any $p\in\mathbb{N}$, lemma \ref{5lemmereference} tells us that
$((\chi_{p}q^{n},\chi_{p}u^{n}))_{n\in\mathbb{N}}$ is
uniformly equicontinuous in
$C(\R^{+};B^{\frac{N}{2}-1}\times(B^{\frac{N}{2}-\frac{3}{2}})^{N})$.
In using Ascoli's theorem we just need to show that
$((\chi_{p}q^{n}(t,\cdot),\chi_{p}u^{n})(t,\cdot))_{n\in\mathbb{N}}$
is relatively compact
in $B^{\frac{N}{2}-1}\times(B^{\frac{N}{2}-\frac{3}{2}})^{N}$ $\forall t\in[0,p\,]$.\\
Let us observe now that the application $u\rightarrow\chi_{p}u$ is
compact from $\widetilde{B}^{\frac{N}{2}-1,\frac{N}{2}}=B^{\N}\cap
B^{\N-1}$ into $\dot{H}^{\frac{N}{2}-1}$, and from
$B^{\frac{N}{2}-1}\cap B^{\frac{N}{2}-\frac{3}{2}}$ into
$\dot{H}^{\frac{N}{2}-\frac{3}{2}}$.\\
After we apply Ascoli's theorem to the family
$((\chi_{p}q^{n},\chi_{p}u^{n}))_{n\in\mathbb{N}}$  on
the time interval $[0,p]$. We then use Cantor's diagonal
process.This finally provides us with a distribution
$(q,u)$ belonging to
$C(\R^{+};\dot{H}^{\frac{N}{2}-1}\times(\dot{H}^{\frac{N}{2}-\frac{3}{2}})^{N})$
and a subsequence (which we still denote by
$(q^{n},u^{n})_{n\in\mathbb{N}}$ such that, for all
$p\in\mathbb{N}$, we have:
\begin{equation}
(\chi_{p}q^{n},\chi_{p}u^{n})\rightarrow_{n\mapsto+\infty}(\chi_{p}q,\chi_{p}u)
\;\mbox{in}\;C([0,p];\,\dot{H}^{\frac{N}{2}-1}\times(\dot{H}^{\frac{N}{2}-\frac{3}{2}})^{N})
\label{5converg}
\end{equation}
This obviously entails that $(q^{n},u^{n})$ tends to
$(q,u)$ in ${\cal D}^{'}(\R^{+}\times\R^{N})$.\\
\\
Coming back to the uniform estimates of step one, we moreover get
that $(q,u)$ belongs to:\\
$$L^{1}(\widetilde{B}^{\frac{N}{2}-1,\frac{N}{2}}\times(B^{\frac{N}{2}+1})^{N})\cap
L^{\infty}(\widetilde{B}^{\frac{N}{2}-1,\frac{N}{2}}\times(B^{\frac{N}{2}+1})^{N})$$
and to
$C^{\frac{1}{2}}(\R^{+};B^{\frac{N}{2}-1})\times(C^{\frac{1}{4}}(\R^{+};B^{\frac{N}{2}-\frac{3}{2}})^{N})$.
Obviously, we have the bounds provided of the firts step.
\\
Let us now prove that $(q,u)$ solves
the system $(SW)$, we first recall that $(q^{n},u^{n})$ solves the
following system:
$$
\begin{cases}
\begin{aligned}
&\p_{t}q^{n}+J_{n}(u^{n}\cdot\n q^{n})+{\rm div}u^{n}=-J_{n}(q^{n}{\rm
div}u^{n})\\
&\p_{t}u^{n}-\bar{\nu}\D u^{n}+\del\n q^{n}-\ka\phi*\n
q^{n}+J_{n}(u^{n}\cdot\n
u^{n})+J_{n}(K(\va(\bar{\rho}(1+q^{n}))\n q^{n})\\
&\hspace{9cm}+J_{n}({\cal A}(\va(\bar{\rho}(1+q^{n})),u^{n}))=0\\
\end{aligned}
\end{cases}
$$
The only problem is to pass to the limit in
${\cal{D}}^{'}(\R^{+}\times\R^{N})$ in the non linear terms. This can
be done by using the convergence results coming from the uniform
estimates (\ref{5converg}).\\
As it is just a matter of doing tedious verifications, we show as a
example the case
of the term $J_{n}(K(\va(\bar{\rho}(1+q^{n})))\n q^{n})$ and $J_{n}({\cal{A}}(\va(\bar{\rho}(1+q^{n})),u^{n}))$.\\
We decompose:
$$J_{n}(K(\va(\bar{\rho}(1+q^{n})))\n q^{n})-K(\rho^{n})\n q^{n}=J_{n}(K(\va(\bar{\rho}(1+q^{n})))\n q^{n})-K(\va(\bar{\rho}(1+q)))\n q.$$
(Note that  for $n$ big enough, we have $K(\va(\bar{\rho}(1+q^{n})))=K(\rho^{n})$ as we control
$\|\rho^{n}\|_{L^{\infty}}$ and $\|\frac{1}{\rho^{n}}\|_{L^{\infty}}$). Next we have:
$$
\begin{aligned}
&J_{n}(K(\va(\bar{\rho}(1+q^{n})))\n q^{n})-K(\va(\bar{\rho}(1+q)))\n q=J_{n}A_{n}+(J_{n}-I)K(\va(\bar{\rho}(1+q)))\n q,\\[2mm]
&\hspace{4cm}\mbox{where}\;\;\;A_{n}=K(\va(\bar{\rho}(1+q^{n})))\n q^{n}-
K(\va(\bar{\rho}(1+q)))\n q.\\
\end{aligned}
$$
We have then $(J_{n}-I)K(\va(\bar{\rho}(1+q)))\n q$ tends to zero as $n\rightarrow+\infty$ due to the property of $J_{n}$ and the fact that
$K(\va(\bar{\rho}(1+q)))\n q$ belongs to $L^{\infty}(B^{\N-1})\hookrightarrow L^{\infty}(L^{q})$ for some $q\geq2$.
Choose $\psi\in C_{0}^{\infty}([0,T)\times\R^{N})$ and $\va^{'}\in C_{0}^{\infty}([0,T)\times\R^{N})$ such that $\va^{'}=1$ on $\mbox{supp}\,\psi$, we
have:
$$|<(J_{n}-I)K(\va(\bar{\rho}(1+q)))\n q,\psi>|\leq\|\va^{'}\, K(\va(\bar{\rho}(1+q)))\n q\|_{L^{\infty}(L^{2})}\|(J_{n}-I)\psi\|_{L^{2}},$$
because $L^{q}_{loc}\hookrightarrow L^{2}_{loc}$ and we conclude by the fact that $\|(J_{n}-I)\psi\|_{L^{2}}\rightarrow0$ as $n$ tends to $+\infty$.\\
Next:
$$
\begin{aligned}
<J_{n}A_{n},\psi>=I^{1}_{n}+I^{2}_{n},
\end{aligned}
$$
with:
$$
\begin{aligned}
&I^{1}_{n}=<(K(\va(\bar{\rho}(1+q^{n})))-K(\va(\bar{\rho}(1+q))))\n q^{n},J_{n}\psi>,\\
&I^{2}_{n}=<K(\va(\bar{\rho}(1+q)))\n(q^{n}-q),J_{n}\psi>.\\
\end{aligned}
$$
We have then:
$$I^{1}_{n}\leq\|\va^{'}q^{n}\|_{L^{\infty}(B^{\N})}\|\va^{'}(q^{n}-q)\|_{L^{\infty}(\dot{H}^{\N-1})}\|\psi\|_{L^{\infty}},$$
Indeed we just use the fact that $\va^{'} B^{\N-1}$ and $\va^{'}\dot{H}^{\N-1}$ are embedded in $L^{2}$. Next we conclude as
we have seen that $q^{n}\rightarrow_{n\rightarrow+\infty} q$ in $C_{loc}(H_{loc}^{\N-1})$. So we obtain:
$$I^{1}_{n}\rightarrow_{n\rightarrow+\infty}0\;\;\;\mbox{in}\;\;{\cal D}^{'}((0,T^{*})\times\R^{N}).$$
We proceed similarly for $I^{2}_{n}$, indeed we have:
$$I^{2}_{n}=<\va^{'} (q^{n}-q),\va^{'}{\rm div}(K(\va(\bar{\rho}(1+q)))J_{n}\psi)>$$
and we have
$K(\va(\bar{\rho}(1+q)))J_{n}\psi\in L^{\infty}(B^{\N})$ so:
$$I^{2}_{n}\leq\|\va^{'} (q^{n}-q)\|_{L^{\infty}(\dot{H}^{\N-1})}\|K(\va(\bar{\rho}(1+q)))J_{n}\psi\|_{L^{\infty}(B^{\N})}.$$
We conclude then that:
$$I^{2}_{n}\rightarrow_{n\rightarrow+\infty}0\;\;\;\mbox{in}\;\;{\cal D}^{'}((0,T^{*})\times\R^{N}).$$
We concentrate us now on the term $J_{n}({\cal{A}}(\va(\bar{\rho}(1+q^{n})),u^{n}))$. Let $\va^{'}\in C_{0}^{\infty}(\R^{+}\times\R^{N})$ and
$p\in\mathbb{N}$
be such that ${\rm supp}\va^{'}\subset[0,p]\times B(0,p)$.
We use the decomposition for $n$ big enough:
$$
\begin{aligned}
\va^{'}J_{n}{\cal A}(\va(\bar{\rho}(1+q^{n})),u^{n})-\va^{'}{\cal A}(\rho,u)=\va^{'}
\chi_{p}{\cal{A}}&(\va(\bar{\rho}(1+q^{n})),\chi_{p}(u^{n}-u))\\
&+\va^{'}{\cal{A}}(\chi_{p}\va(\bar{\rho}(1+q^{n}))
-\chi_{p}\bar{\rho}(1+q)),u).\\
\end{aligned}
$$
According to the uniform
estimates and (\ref{5converg}), $\chi_{p}(u^{n}-u)$
tends to $0$ in $L^{1}([0,p];\dot{H}^{\N+1})$  by interpolation so that the first
term tends to $0$ in $L^{1}(\dot{H}^{\N-1})$ and we conclude for the second term
in remarking that $\frac{\va}{\rho_{n}}$ tends to $\frac{\va}{\rho}$ as $\rho_{n}$
in $L^{\infty}(L^{\infty}\cap \dot{H}^{\N})$.\\
The other nonlinear terms can be treated in the same way.
\subsection{Proof of the uniqueness in the critical case}
\begin{theorem}
\label{unicitecritique}
Let $N\geq2$, and $(q_{1},u_{1})$ and $(q_{2},u_{2})$ be solutions of $(SW)$ with the same data $(q_{0},u_{0})$
on the time interval $[0,T^{*})$. Assume that for $i=1,2$:
$$(q_{i},u_{i})\in C([0,T^{*}),B^{1}_{N,1})\;\;\;\mbox{and}\;\;\;u_{i}\in \big(C([0,T^{*}),B^{0}_{N,1})\cap
L^{1}_{loc}([0,T^{*}),B^{2}_{N,1})\big)^{N}.$$
There exists a constant $\alpha>0$ depending only on $N$ and physical constants such that if:
\begin{equation}
\|q_{1}\|_{\widetilde{L}^{\infty}_{T^{*}}(B^{1}_{N,1})}\leq\alpha,
\label{5estimimsurq}
\end{equation}
then $(q_{1},u_{1})=(q_{2},u_{2})$ on $[0,T^{*})$.
\end{theorem}
Let $(q_{1},u_{1})$, $(q_{2},u_{2})$ belong to $E^{\N}$ with the
same initial data, we set $(\delta q,\delta
u)=(q_{2}-q_{1},u_{2}-u_{1})$. We can then write the system $(SW)$
as follows:
\begin{equation}
\begin{cases}
\begin{aligned}
&\frac{\p}{\p t}\delta q+u_{2}\cdot\n\delta q=H_{1},\\[2mm]
&\frac{\p}{\p t}\delta u-\bar{\nu}\D\delta u=H_{2}\\
\end{aligned}
\end{cases}
\label{5systemunicite}
\end{equation}
with:
$$
\begin{aligned}
&H_{1}=-{\rm div}\delta u-\delta u\cdot\n q_{1}-\delta q{\rm div}u_{2}-q_{1}{\rm div}u,\\[2mm]
&H_{2}=-\del\n \delta q-\ka \phi*\n\delta q-u_{2}\cdot\n\delta u-\delta
u\cdot\n u_{1}+{\cal A}(q_{1},\delta u)+
{\cal A}(\delta q,u_{2}).\\
\end{aligned}
$$
Due to the term $\de u\cdot\n q^{1}$ in the right-hand side of the
first equation, we loose one derivative when estimating $\de q$:
one only gets bounds in $L^{\infty}(B^{0}_{N,1})$.\\
Now, the right hand-side of the second equation contains a term
of type ${\cal A}(\delta q,u_{2})$ so that the loss of one
derivative for $\de q$ entails a loss of one derivative for $\de u$.
Therefore, getting bounds in:
$$C(\R^{+};B^{-1}_{N,1})\cap L^{1}(\R^{+};B^{1}_{N,1})$$
for $\de u$ is the best that one can hope. If enough regularity were
available, we would not have to worry about this loss of derivative.
But in the present case, the above heuristic fails because we have
reached some limit cases for the product laws. Indeed, a term such
as $\de u\cdot\n u_{1}$ cannot be estimated properly: the product
does not map $B^{0}_{N,1}\times B^{0}_{N,1}$ into $ B^{-1}_{N,1}$
but in the somewhat {\it larger space} $ B^{-1}_{N,\infty}$. At this
point, we could try instead to get bounds for $\de u$ in:
$$C([0,T^{*});B^{-1}_{N,\infty})\cap L^{1}_{loc}([0,T^{*});B^{1}_{N,\infty}),$$
but we then have to face the lack of control on $\de u$ in
$L^{1}(0,T;L^{\infty})$ (because in contrast with $B^{1}_{N,1}$, the
space $B^{1}_{N,\infty}$ is not imbedded in $L^{\infty}$) so that we
run into troubles when estimating $\de u\cdot\n q_{1}$. The key to
that difficulty relies on the following logarithmic interpolation
inequality (see the proposition \ref{5interp1}):
$$\|u\|_{L^{1}_{T}(B^{1}_{N,1})}\lesssim\|u\|_{\widetilde{L}^{1}_{T}(B^{1}_{N,\infty})}
\log\biggl(e+\frac{\|u\|_{\widetilde{L}^{1}_{T}(B^{0}_{N,\infty})}+\|u\|_{\widetilde{L}^{1}_{T}(B^{2}_{N,\infty})}}
{\|u\|_{\widetilde{L}^{1}_{T}(B^{1}_{N,\infty})}}\biggl),
$$
and a well-known generalization of Gr\"onwall the Osgood's lemma
(see \cite{5DFourier}) that we recall.
\begin{lemme}
\label{5osgood}
Let $F$ be a measurable positive function and $\gamma$ a positive
locally integrable function, each defined on the domain
$[t_{0},t_{1}]$. Let $\mu:\;[0,+\infty)\rightarrow[0,+\infty)$ be a
continuous nondecreasing function, with $\mu(0)=0$. Let $a\geq0$,
and assume that for all $t\in[t_{0},t_{1}]$,
$$F(t)\leq a+\int^{t}_{t_{0}}\gamma(s)\mu(F(s))ds.$$
If $a>0$, then:
$$-{\cal M}(F(t))+{\cal M}(a)\leq\int^{t}_{t_{0}}\gamma(s)ds,\;\;\;\mbox{where}\;\;
{\cal M}(x)=\int^{1}_{x}\frac{ds}{\mu(s)}.$$ If $a=0$ and ${\cal
M}(0)=+\infty$, then $F=0$.
\end{lemme}
{\bf Proof of the theorem \ref{unicitecritique}:}
\subsubsection*{First step: in which space do we work?}
Let us observe first that in view of Sobolev embedding, $q_{i}\in
C(\R^{+};L^{\infty})$. Therefore, if $\alpha$ is small enough,
by embedding and continuity we get:
$$|q_{i}(t,x)|\leq\frac{1}{2}$$
for $x\in\R^{N}$ and $t$ in a small nontrivial time interval
$[0,T]$.\\
That observation will enable us to apply proposition \ref{5composition} to the
non-linear terms involving $q_{i}$. \\
We shall further assume that
$T\in (0,+\infty)$ has been chosen so small as to satisfy:
\begin{equation}
C\|\n u_{2}\|_{L^{1}_{T}(B^{1}_{N,1})}\leq\log2,
\label{5controlu}
\end{equation}
for some appropriate constant $C$ whose meaning will be clear from
the computations below.\\
To begin with, we shall prove uniqueness on the time interval
$[0,T]$ by estimating $(\de q,\de u)$ in the following functional
space:
$$F_{T}=L^{\infty}([0,T];B^{0}_{N,\infty})\times(L^{\infty}([0,T];B^{-1}_{N,\infty})\cap
\widetilde{L}^{1}_{T}(B^{1}_{N,\infty}))^{N}.$$
Indeed as explained below, in this space we can control the remainder because it is appropriate to the result of paraproduct.
\subsubsection*{Why $(\de q,\de u)$ is in $F_{T}$?}
Of course, we have
to state that $(\de q,\de u)\in F_{T}$, a fact which is not entirely
obvious. We want now to show that $(\de q,\de u)$ belongs to $F_{T}$.\\
According to our assumption on $(q_{i},u_{i})$, the estimates of paraproduct
yield $\p_{t}q_{i}\in L^{2}_{T}(B^{0}_{N,1})$. Therefore
$\bar{q}_{i}=q_{i}-q_{0}$ belongs to
$C^{\frac{1}{2}}([0,T],B^{0}_{N,1})$, which clearly entails by embedding $\de
q\in C([0,T],B^{0}_{N,\infty})$.\\
Let $\bar{u}_{i}=u_{i}-u_{L}$ with
$u_{L}$ solution to the following linear heat equation:
$$
\begin{cases}
\begin{aligned}
&\p_{t}u_{L}-\mu \D u_{L}=-\del\n q_{0}+\ka\n(\phi*q_{0}),\\
&u_{L}(0)=u_{0}.\\
\end{aligned}
\end{cases}
$$
We obviously have $(\bar{u}_{i})_{0}=0$ and:
$$\p_{t}\bar{u}_{i}-\bar{\nu}\D\bar{u}_{i}=-\bar{u}_{i}\cdot\n\bar{u}_{i}-\del\n q_{i}+\ka\n(\phi*q_{i})
+{\cal A}(\rho_{i},q_{i})+K(\rho_{i},u_{i}).$$ The product and
composition laws in Besov spaces insure that the right-hand side
belongs to $L^{2}_{T}(B^{-1}_{N,\infty})$ (because
$B^{0}_{N,1}\times B^{0}_{N,\infty}\hookrightarrow B^{-1}_{N,\infty}$) thus to
$\widetilde{L}^{1}_{T}(B^{-1}_{N,\infty}))$ (for the last term, we
use that $\bar{q}_{i}\in L^{\infty}_{T}(B^{0}_{N,1})$.\\
Now Proposition \ref{5chaleur} implies that:
$$\bar{u}_{i}\in L^{\infty}_{T}(B^{-1}_{N,\infty})\cap
\widetilde{L}^{1}_{T}(B^{1}_{N,\infty})).$$
\subsubsection*{Second step: Estimates on $(\de q,\de u)$}
Let us turn to estimate $\de q$. Proposition \ref{5transport} combined with (\ref{5controlu}) yields for $t\leq T$:
$$\|\de q\|_{L^{\infty}_{t}(B^{0}_{N,\infty})}\lesssim
\int^{t}_{0}\big(\|\de u\cdot\n q\|_{B^{0}_{N,\infty}}+\|\de q\;{\rm
div}u_{2}\|_{B^{0}_{N,\infty}}+\|{\rm div}\de
u\|_{B^{0}_{N,\infty}})\big)d\tau
$$
Estimate of type $B^{\NN}_{p,\infty}\cap L^{\infty}\times B^{s}_{p,\infty}\hookrightarrow B^{s}_{p,\infty}$ with $s+\NN>0$ enables us to get the
following inequality:
$$\|\de q\|_{L^{\infty}_{t}(B^{0}_{N,\infty})}\lesssim
\int^{t}_{0}\big(\|\de q\|_{B^{0}_{N,\infty}}\|{\rm
div}u_{2}\|_{B^{1}_{N,\infty}\cap L^{\infty}}+\|\de
u\|_{B^{1}_{N,\infty}\cap L^{\infty}}(1+\|
q_{1}\|_{B^{1}_{N,1}}\big)d\tau,$$
whence, according to Gronwall
inequality, to the embedding $B^{1}_{N,1}\hookrightarrow
B^{1}_{N,\infty}\cap L^{\infty}$ and to (\ref{5estimimsurq}) we get:
$$\|\de q\|_{L^{\infty}_{t}(B^{0}_{N,\infty})}\lesssim\|\de u\|_{L^{1}_{t}(B^{1}_{N,1})}
(1+\|q_{1}\|_{L^{\infty}_{t}(B^{1}_{N,1})}).$$
Making use of (\ref{5estimimsurq}) and proposition \ref{5interp}, we end up with:
$$\|\de q\|_{L^{1}_{t}(B^{0}_{N,\infty})}\lesssim\|\de u\|_{\widetilde{L}^{1}_{t}(B^{1}_{N,\infty})}
\log\biggl(e+\frac{\|\de
u\|_{\widetilde{L}^{1}_{t}(B^{0}_{N,\infty})}+\|\de
u\|_{\widetilde{L}^{1}_{t}(B^{2}_{N,\infty})}} {\|\de
u\|_{\widetilde{L}^{1}_{t}(B^{1}_{N,\infty})}}\biggl).
$$
Remark that:
$$\|\de
u\|_{\widetilde{L}^{1}_{t}(B^{0}_{N,\infty})}+\|\de
u\|_{\widetilde{L}^{1}_{t}(B^{2}_{N,\infty})}\leq
V(t)=V_{1}(t)+V_{2}(t)$$ with:
$$V_{i}(t)=\int^{t}_{0}(\|u_{i}(\tau)\|_{B^{0}_{N,1}}+\|u_{i}(\tau)\|_{B^{2}_{N,1}})d\tau<+\infty$$
since
$\widetilde{L}^{\infty}_{t}(B^{0}_{N,1})\hookrightarrow\widetilde{L}^{1}_{t}(B^{0}_{N,1})$
for finite $t$.\\
We finally get:
\begin{equation}
\|\de q\|_{L^{1}_{t}(B^{0}_{N,\infty})}\lesssim\|\de u\|_{\widetilde{L}^{1}_{t}(B^{1}_{N,\infty})}
\log\biggl(e+\frac{V(t)} {\|\de
u\|_{\widetilde{L}^{1}_{t}(B^{1}_{N,\infty})}}\biggl),
\label{5controldev}
\end{equation}
with $V$
non-decreasing bounded function of $t\in[0,+\infty)$.\\
\\
Let us now turn to the proof of estimates for $\de u$. According to
proposition \ref{5chaleur}, we have:
$$
\begin{aligned}
&\|\de u\|_{L^{\infty}_{t}(B^{-1}_{N,\infty})}+\|\de
u\|_{L^{1}_{t}(B^{1}_{N,\infty})}\lesssim\|u_{2}\cdot\n\de
u\|_{\widetilde{L}^{1}_{t}(B^{-1}_{N,\infty})}+\|{\cal A}(q_{1},\de
u)\|_{\widetilde{L}^{1}_{t}(B^{-1}_{N,\infty})}\\[2mm]
&+\|\de u\cdot\n
u_{1}\|_{\widetilde{L}^{1}_{t}(B^{-1}_{N,\infty})}+\|{\cal A}(\de
q,u_{2})\|_{\widetilde{L}^{1}_{t}(B^{-1}_{N,\infty})}
+\|K(\de q)\n q_{2}\|_{\widetilde{L}^{1}_{t}(B^{-1}_{N,\infty})}\\[2mm]
&\hspace{10cm}+\|K(q_{2})\n\de q\|_{\widetilde{L}^{1}_{t}(B^{-1}_{N,\infty})}.\\
\end{aligned}
$$
Let us assume that the $\alpha$ appearing in (\ref{5estimimsurq}) is small enough so that the second term
 in the right-hand side may be absorbed by the left-hand side. $\|u_{2}\|_{\widetilde{L}^{2}_{t}(B^{1}_{N,1})}$ tends to $0$ when $t$ goes to $0$,
  so if we choose $T$
small enough, the first term may also be absorbed.
Using interpolation of proposition \ref{5interp}, we obtain for all $t\in[0,T]$,
$$
\begin{aligned}
&\|\de u\|_{L^{\infty}_{t}(B^{-1}_{N,\infty})}+\|\de
u\|_{L^{1}_{t}(B^{1}_{N,\infty})}\lesssim\int^{t}_{0}\big[\|
u_{1}\|_{B^{2}_{N,1}}\| \de
u\|_{B^{-1}_{N,\infty}}+(1+\|u_{2}\|_{B^{2}_{N,1}})\|
\de q\|_{B^{0}_{N,\infty}}\big]d\tau\\
\end{aligned}
$$
Let us now plug (\ref{5controldev}) in the above inequality. Denoting:
$$X(t)=\|\de u\|_{L^{\infty}_{t}(B^{-1}_{N,\infty})}+\|\de
u\|_{L^{1}_{t}(B^{1}_{N,\infty})}$$ we get for $t\leq T$,
$$
\begin{aligned}
X(t)&\lesssim\int^{t}_{0}\big(1+V^{'}(\tau)\big)X(\tau)\log\biggl(e+\frac{V(\tau)}{X(\tau)}\biggl)d\tau\\
\end{aligned}
$$
As:
$$V^{'}\in L^{1}(0,T)\;\;\;\mbox{and}\;\;\;\int^{1}_{0}\frac{dr}{r\log(e+\frac{V(T)}{r}}=+\infty,$$
Osgood's lemma (see lemma \ref{5osgood}) implies that $X=0$ on $[0,T]$, whence also $\de
q=0$. Standard arguments of connexity then yield uniqueness on the whole interval
$[0,+\infty)$.
\section{Proof of theorem \ref{5theo2}}
\label{5section5} We will proceed similarly to the proof of the
theorem \ref{5theo1}. To begin with, let us observe that under the
definition \ref{5bar}, system reads:
\begin{equation}
\begin{cases}
\begin{aligned}
&\p_{t}q+u\cdot\n u=H\\
&\p_{t}u-\bar\nu \D u=-K\\
&(q,u)_{t=0}=(q_{0},u_{0})\\
\end{aligned}
\end{cases}
\label{5syst4}
\end{equation}
with:
$$
\begin{aligned}
&H=-(1+q){\rm div}u,\\
&K=G-u\cdot\n u+(P^{'}(\bar{\rho})+\kappa)\n\rho-\kappa\phi*\n\rho.\\
\end{aligned}
$$
As previously we can build approximate smooth solutions
$(q^{n},u^{n})$ of (\ref{5syst4}) in studying the Korteweg system with a capillarity
coefficient $\kappa_{n}=\frac{1}{n}$. It is convenient to split $(q^{n},u^{n})$
into the solution of a linear system with initial data $(q_{0}^{n},u_{0}^{n})$, and the discrepancy to that
solution. More precisely we denote by $(q^{n}_{L},u^{n}_{L})$ the solution of the linearized pressure-less system on the intervall $[0,T]$:
$$
\begin{cases}
\begin{aligned}
&\p_{t}q^{n}_{L}+{\rm div}u^{n}_{L}=0,\\
&\p_{t}u^{n}_{L}-\bar{\mu}\D u^{n}_{L}-(\bar{\lambda}+\bar{\mu})\n{\rm div}u^{n}_{L}=0,\\
&(q^{n}_{L},u^{n}_{L})_{/t=0}=(q_{0}^{n},u_{0}^{n}),\\
\end{aligned}
\end{cases}
$$
with:
$$(q_{0}^{n},u_{0}^{n})=(\sum_{|l|\leq n}\D_{l}q_{0},\sum_{|l|\leq n}\D_{l}u_{0}).$$
We set:
$$(q^{n},u^{n})=(q^{n}_{L}+\bar{q}^{n},u^{n}_{L}+\bar{u}^{n})$$
We can state now that $(q^{n},\bar{u}^{n})$ verifies the following
linear system:
\begin{equation}
\begin{cases}
\begin{aligned}
&\p_{t}\bar{q}^{n}+{\rm div}\bar{u}^{n}=F^{n},\\
&\p_{t}\bar{u}^{n}-\bar{\mu}\D \bar{u}^{n}-(\bar{\lambda}+\bar{\mu})\n{\rm div}\bar{u}^{n}-\frac{1}{n}\n\D \bar{q}^{n}=G^{n},\\
&(\bar{q}^{n},\bar{u}^{n})_{/t=0}=(0,0)\\
\end{aligned}
\end{cases}
\label{5syst4}
\end{equation}
with:
$$
\begin{aligned}
&F^{n}=-q^{n}{\rm div}u^{n},\\
&G^{n}=-u^{n}\cdot\n u^{n}+{\cal A}(\rho^{n},u^{n})-K(\rho^{n})\n
q^{n}+ (P^{'}(\bar{\rho})+\kappa)\n\rho^{n}-\kappa\phi*\n\rho^{n}.
\end{aligned}
$$
We want show that such solution $(q^{n},u^{n})$ exists, in this goal we recall some theorem by R. Danchin and B. Desjardins in \cite{5DD}.
\begin{theorem}
Let $p\in[1,+\infty[$. Then there exists $\eta>0$ such that if $q_{0}\in B^{\NN}$, $u_{0}\in(B^{\NN-1})^{N}$ and:
$$\|q_{0}\|_{B^{\NN}}\leq\eta,$$
then there exists $T>0$ such that system (\ref{5syst4}) has a unique solution $(q,u)$ in $\widetilde{E}^{p,\kappa}_{T}$.
\end{theorem}
In fact, we can extend this result to the case where the viscosity coefficients are variable. And we can show that
the uniform estimates are independent of the capillarity coefficient. We can obtain the following result on our solution
$(q^{n},u^{n})_{n\in\mathbb{N}}$.
\begin{theorem}
\label{5theoutile}
Let $p\in[1,+\infty[$. Then there exists $\eta>0$ such that if $q^{n}_{0}\in B^{\NN}$, $u^{n}_{0}\in(B^{\NN-1})^{N}$ and:
$$\|q^{n}_{0}\|_{B^{\NN}}\leq\eta,$$
then there exists $T>0$ such that system (\ref{5syst4}) has a unique solution $(q^{n},u^{n})$ in $\widetilde{E}^{p,\frac{1}{n}}_{T}$ and
$(q^{n},u^{n})$ are uniformly bounded in $F^{p}_{T}$.
\end{theorem}
{\bf Proof:}
\\
\\
We recall a proposition (see \cite{5DD}) on the following linearized pressure-less system:
\begin{equation}
\begin{aligned}
&\p_{t}q+{\rm div}u=F,\\
&\p_{t}u-\bar{\mu}\D u-(\bar{\lambda}+\bar{\mu})\n{\rm div}u-\bar{\kappa}\n\D q=G.
\end{aligned}
\label{5systemkorteweg}
\end{equation}
\begin{proposition}
\label{5proprappel}
Let $s\in\R$, $p\in[1,+\infty]$, $1\leq\rho_{1}\leq+\infty$ and $T\in]0,+\infty]$. If $(q_{0},u_{0})\in B^{s}_{p}\times (B^{s-1}_{p})^{N}$ and
$(F,G)\in\widetilde{L}_{T}^{\rho_{1}}(B^{s-2+\frac{2}{\rho_{1}}}_{p}\times
(B^{s-3+\frac{2}{\rho_{1}}}_{p})^{N})$, then the  above linear system has a unique solution
$(q,u)\in\widetilde{C}_{T}(B^{s}_{p}\times(B^{s-1}_{p})^{N})\cap\widetilde{L}_{T}^{\rho_{1}}(B^{s+\frac{2}{\rho_{1}}}_{p}\times
(B^{s-1+\frac{2}{\rho_{1}}}_{p})^{N})$.\\
Moreover for all $\rho\in[\rho_{1},+\infty]$, there exists a constant $C$ depending only on $\bar{\mu}$, $\bar{\lambda}$, $\bar{\kappa}$, $p$,
$\rho_{1}$ and $N$ such that the following inequality holds:
$$
\begin{aligned}
\|q\|_{\widetilde{L}_{T}^{\rho}(B^{s+\frac{2}{\rho}}_{p})}+\|u\|_{\widetilde{L}_{T}^{\rho}(B^{s-1+\frac{2}{\rho})}_{p}}\leq
C(\|q_{0}&\|_{B^{s}_{p}}+\|u_{0}\|_{B^{s-1}_{p}}\\
&+\|F\|_{\widetilde{L}_{T}^{\rho_{1}}(B^{s-2+\frac{2}{\rho_{1}}}_{p})}
+\|G\|_{\widetilde{L}_{T}^{\rho_{1}}(B^{s-3+\frac{2}{\rho_{1}}}_{p})}).\\
\end{aligned}
$$
\end{proposition}
\begin{remarka}
More precisely we have:
$$
\begin{aligned}
&\|q\|_{\widetilde{L}_{T}^{\infty}(B^{s}_{p})}+\kappa\|q\|_{\widetilde{L}_{T}^{1}(B^{s+2}_{p})}\leq
C(\|q_{0}\|_{B^{s}_{p}}+\|u_{0}\|_{B^{s-1}_{p}}+\|F\|_{\widetilde{L}_{T}^{1}(B^{s}_{p})}
+\|G\|_{\widetilde{L}_{T}^{1}(B^{s-1}_{p})}).
\end{aligned}
$$
\end{remarka}
\subsubsection*{Uniform Estimates for $(q^{n},u^{n})_{n\in\mathbb{N}}$}
The existence of solutions in the case of general viscosity coefficient follow the same line as the proof in \cite{5DD}. It suffices to solve a problem
of fixed point.\\
Denoting by $V(t)$ the semi-group generated by system (\ref{5systemkorteweg}), we have:
$$(q^{n}_{L},u^{n}_{L})(t)=V(t)(q^{n}_{0},u^{n}_{0}).$$
Let us define:
$$\phi_{q^{n}_{L},u^{n}_{L}}(\bar{q}^{n},\bar{u}^{n})=\int^{t}_{0}V(t-s)(F(q^{n}_{L}+\bar{q}^{n},u_{L}^{n}+\bar{u}^{n})(s),G(q^{n}_{L}+\bar{q}^{n},u_{L}^{n}+\bar{u}^{n})(s))\,ds.$$
where we have set:
$$
\begin{aligned}
&F(q,u)=-{\rm div}(qu),\\
&G(q,u)=-u\cdot\n u++{\cal A}(\rho,u)-K(\rho)\n
q+ (P^{'}(\bar{\rho})+\kappa)\n\rho-\kappa\phi*\n\rho.\\
\end{aligned}
$$
In order to prove the existence part of the theorem \ref{5theoutile}, ther's just have to show that $\phi_{q^{n}_{L},u^{n}_{L}}$
has a fixed point in $F_{T}^{p}$. Since $F_{T}^{p}$ is a Banach space, we can prove that $\phi_{q^{n}_{L},u^{n}_{L}}$ satisfies the
 hypothesis of Picard's theorem in a ball $(B(0,R)$ of $F_{T}^{p}$ for sufficiently small $R$. Moreover $R$ depends only of mathematical constants. We
 can find a time $T$ independent of $n$ such that for all initial data verifying $\|q^{n}_{0}\|_{B^{\NN}}\leq\frac{R}{2}$, we have existence of
 solution $(q^{n},u^{n}$ at less on the interval $(0,T)$. The end of proof consists to verify that $(q^{n},u^{n})$ is uniform in $F^{p}_{T}$, it
 suffices to use the proposition \ref{5proprappel}.\\
We now have builded approximated solution $(q^{n},u^{n})$ of system $(SW)$ and we can conclude in using technics of compactness.
\subsubsection*{2) Existence of a solution}
The existence of a solution stems from compactness properties for
the sequence $(q^{n},u^{n})_{n\in\mathbb{N}}$ an we want use some
result of type Ascoli as in the proof of theorem \ref{5theo1}.
\begin{lemme}
\label{5lemmeutile} The sequence $(\p_{t}\bar{q}^{n},\p_{t}\bar{u}^{n})_{n\in\mathbb{N}}$ is
uniformly bounded in:
$$L^{2}(0,T;\widetilde{B}^{\NN,\NN-1})\times(L^{\alpha}(0,T;\widetilde{B}^{\NN-1,\NN-2}))^{N},$$
for some $\alpha>1$.
\end{lemme}
{\bf Proof:}
\\
\\
Throughout the proof, we will extensively use that
$\widetilde{L}^{\rho}_{T}(B^{s}_{p})\hookrightarrow
L^{\rho}_{T}(B^{s}_{p})$. The notation u.b will stand for uniformly bounded.\\
We have:
\begin{equation}
\begin{aligned}
&\p_{t}q^{n}=-u^{n}\cdot\n q^{n}-(1+q^{n}){\rm div}u^{n},\\
&\p_{t}\bar{u}^{n}=-u^{n}\cdot\n
u^{n}-q^{n}{\cal A}(\rho^{n},\bar{u}^{n})
-K(q^{n})\n q^{n}+\frac{1}{n}\n\D\bar{q}^{n}.\\
\end{aligned}
\label{5T8}
\end{equation}
We start with show that $\p_{t}\bar{q}^{n}$ is u.b in $L^{2}(0,T;\widetilde{B}^{\NN,\NN-1}_{p})$.\\
Since $u^{n}$ is u.b in $L^{2}_{T}(B^{\NN})$ and $\n q^{n}$ is u.b
in $L^{\infty}_{T}(B^{\NN-1})$, then $u^{n}\cdot\n q^{n}$ is u.b
in $L^{2}_{T}(\widetilde{B}^{\NN,\NN-1})$. 
Similar arguments enable us to conclude for the term
$(1+q^{n}){\rm div}u^{n}$ which is u.b in
$L^{2}_{T}(\widetilde{B}^{\NN,\NN-1}_{p})$ because
$q^{n}$ is u.b in $L^{\infty}_{T}(B^{\NN}_{p})$ and ${\rm div}u^{n}$ is u.b in  $L^{2}_{T}(B^{\NN-1}_{p})$.\\
\\
Let us now study $\p_{t}\bar{u}^{n+1}$. According to step one and to
the definition of $u^{n}_{L}$, the term ${\cal A}\bar{u}^{n+1}$ is
u.b in $L^{2}(B^{\NN-2}_{p})$. Since $u^{n}$ is u.b in
$L^{\infty}(B^{\NN-1}_{p})$ and $\n u^{n}$ is u.b in $L^{2}(B^{\NN-1}_{p})$,
so $u^{n}\cdot\n u^{n}$ is u.b in $L^{2}(B^{\NN-2}_{p})$ thus
in $L^{2}(\widetilde{B}^{\NN-1,\NN-2}_{p})$.\\
Moreover we have $q^{n}$ is u.b in $L^{\infty}(B^{\NN}_{p})$ and $q^{n}$
is u.b in $L^{\infty}$, so by proposition \ref{5composition} $\n K_{0}(q^{n})$ is u.b in
$L^{\infty}(B^{\NN-1}_{p})$ thus in
$L^{2}(\widetilde{B}^{\NN-1,\NN-2}_{p})$. This concludes the lemma.
{\hfill $\Box$}
\\
\\
Now, let us turn to the proof of the existence of a solution for the
system $(SW)$. We want now use some results of type Ascoli to
conclude in use the properties of compactness of the
lemma \ref{5lemmeutile}.\\
According lemma \ref{5lemmeutile}, $(q^{n},u^{n})_{n\in\mathbb{N}}$ is
u.b in:
$$C^{\frac{1}{2}}([0,T];\widetilde{B}^{\NN,\NN-1}_{p})
\times(C^{1-\frac{1}{\alpha}}([0,T];\widetilde{B}^{\NN-1,\NN-2}_{p}))^{N},$$
thus is uniformly equicontinuous in
$C(([0,T];\widetilde{B}^{\NN,\NN-1}_{p})\times(\widetilde{B}^{\NN-1,\NN-2}_{p})^{N})$.
On the other hand we have the following result of compactness, for
any $\phi\in C_{0}^{\infty}(\R^{N})$, $s\in\R$, $\de>0$ the
application $u\rightarrow\phi u$ is compact from $B^{s}_{p}$ to
$\widetilde{B}^{s,s-\de}_{p}$. Applying Ascoli's theorem, we infer that
up to an extraction $(q^{n},u^{n})_{n\in\mathbb{N}}$ converges in
${\cal D}^{'}([0,T]\times\R^{N})$ to a limit $(\bar{q},\bar{u})$
which belongs to:
$$C^{\frac{1}{2}}([0,T];\widetilde{B}^{\NN,\NN-1}_{p})\times(C^{1-\frac{1}{\alpha}}([0,T];\widetilde{B}^{\NN-1,\NN-2}_{p}))^{N}$$
Let $(q,u)=(\bar{q},\bar{u})+(q_{0},u_{L})$. Using again uniform
estimates of step one and proceeding as, we gather that $(q,u)$
solves $(SW)$ and belongs to:
$$\bar{\rho}+\widetilde{L}^{\infty}_{T}(\widetilde{B}^{\frac{N}{p},\NN-1}_{p})\times\big(
\widetilde{L}^{1}_{T}(B^{\NN+1}_{p})\cap\widetilde{L}^{\infty}_{T}(B^{\NN-1}_{p})\big)^{N}.$$
Applying proposition, we get the continuity results:
$$\rho-\bar{\rho}\in C([0,T],\widetilde{B}^{\frac{N}{p},\NN-1}_{p}),\;\;\;u\in C([0,T],B^{\NN-1}_{p}).$$
\section{Proof of theorem \ref{5theo3}}
\label{5section6} In this section, we consider the case when the
initial density belongs to $\bar{\rho}+\widetilde{B}^{\NN,\NN+\e}$ and
satisfies $0<\bar{\rho}\leq\rho_{0}$. We consider again the study of
a approximate sequence verifying (\ref{5syst4}) and we proceed as
previously.
\subsubsection*{Construction of approximate solutions}
We consider the following system:
\begin{equation}
\begin{aligned}
&\p_{t}q^{n+1}+u^{n}\cdot\n q^{n+1}=q^{n}{\rm div}u^{n},\\
&\p_{t}u^{n+1}+u^{n}\cdot\n
u^{n+1}-{\cal A}(u^{n+1})
+\n q^{n+1}=G^{n},\\
\end{aligned}
\label{5T8}
\end{equation}
with:
$$G^{n}={\cal A}(u^{n})-{\cal A}(\rho^{n},u^{n})+K(\rho^{n})\n q^{n}.$$
\subsubsection*{Uniform estimates for $(q^{n},u^{n})_{n\in\mathbb{N}}$}
Let us show that the sequence $(q^{n},u^{n})_{n\in\mathbb{N}}$ is
uniformly bounded in provided
that $T$ and $\eta$ have been chosen small enough.\\
Let us remark first that according to proposition, there exists a
universal condition $K$ such that for all $n\in\mathbb{N}$, we have:
$$\|u^{n}_{L}\|_{\widetilde{L}^{\infty}_{T}(B^{\NN})}\leq U_{0},\;\;\;\;\;$$
with $U_{0}=K\|u_{0}\|_{B^{\N}}$,
$$\|u^{n}_{L}\|_{\widetilde{L}^{1}_{T}(B^{\NN+1})}\leq K\sum_{q\in\mathbb{Z}}2^{q(\NN-1)}\|\D_{q}u_{0}\|$$
Suppose that $\|q_{0}\|_{B^{\NN}}\leq\eta$ for a small constant
$\eta$ and let $C_{0}=1+\|q_{0}\|_{B^{\NN}}$. According to , we can
choose a positive time $T$ such that the following property holds
for all $n\in\mathbb{N}$:
$$\|u^{n}_{L}\|_{\widetilde{L}^{\infty}_{T}(B^{\NN})}\leq U_{0},\;\;\;\mbox{and}
\;\;\;\|u^{n}_{L}\|_{\widetilde{L}^{r}_{T}(B^{\NN+1})}\leq
\eta^{\frac{2}{r}}U_{0}^{1-\frac{1}{r}}$$ Let us now show by
induction that the following estimates are satisfied:
$$
\begin{aligned}
&\|q^{n}\|_{\widetilde{L}^{\infty}_{T}(B^{\NN})}\leq\sqrt{\eta}\\
&\|\bar{u}^{n}\|_{\widetilde{L}^{\infty}_{T}(B^{\NN-1})}+\|\bar{u}^{n}\|_{\widetilde{L}^{1}_{T}(B^{\NN+1})}
\leq\eta\\
\end{aligned}
$$
From now, we suppose that $\eta\leq(2C_{1})^{-2}$ where $C_{1}$ is
the norm of the injection $B^{\NN}\hookrightarrow L^{\infty}$. This
ensures us that the following inequality is satisfied:
$$\frac{1}{2}\leq 1+q^{n}\leq\frac{3}{2}$$
According to proposition, we have:
$$\|q^{n+1}\|_{\widetilde{L}^{\infty}_{T}(B^{\NN})}\leq \exp^{C\|u^{n}\|_{L^{1}_{T}(B^{\NN+1})}}
\big(\|q_{0}\|_{B^{\NN}}+\|F^{n}\|_{L^{1}_{T}(B^{\NN})}\big)$$
Moreover we have by Proposition:
$$\|u^{n}\|_{\widetilde{L}^{\infty}_{T}(B^{\NN})}\leq \exp^{C\|u^{n}\|_{L^{1}_{T}(B^{\NN+1})}}
\big(\|q_{0}\|_{B^{\NN}}+\|F^{n}\|_{L^{1}_{T}(B^{\NN})}\big)$$ whence
according to proposition:
$$\|F_{n}\|_{L^{1}_{T}(B^{\NN})}\lesssim T^{1-\frac{1}{r_{1}}}(1+\|q^{n}\|_{L^{\infty}_{T}(B^{\NN})})\,
\|{\rm div}u^{n}\|_{L^{r_{1}}_{T}(B^{\NN})}\lesssim
(1+U_{0})^{\frac{1}{2}}\eta,$$ and we get finally:
$$\|q^{n}\|_{\widetilde{L}^{\infty}_{T}(B^{\NN})}\leq (1+U_{0})^{\frac{1}{2}}\eta\exp^{C(1+U_{0})^{\frac{1}{2}}\eta}$$
Obviously if $\eta$ has been chosen small enough then $q^{n+1}$ satisfies the estimate in $({\cal P}_{n+1})$.\\
Applying proposition to the second equation of yields:
$$
\begin{aligned}
&\|\bar{u}^{n+1}\|_{\widetilde{L}^{\infty}_{T}(B^{\NN})}+\|\bar{u}^{n+1}\|_{L^{1}_{T}(B^{\NN+1})}\lesssim
\|u^{n}\cdot\n u^{n}\|_{L^{1}_{T}(B^{\NN-1})}+\|K_{0}(q^{n})\n q^{n}\|_{L^{1}_{T}(B^{\NN-1})}\\[1,8mm]
&\hspace{10cm}+\|{\cal A}(q^{n},u^{n})\|_{L^{1}_{T}(B^{\NN-1})}.\\
\end{aligned}
$$
We now use proposition \ref{5composition} as in the
proof of \ref{5T3} to conclude.
\subsubsection*{Existence of a solution}
We can now easily show that $(q^{n},u^{n})$ is a Cauchy sequel in our space $F_{T}$ of uniqueness and so
$(q_{n},u_{n})\rightarrow(q,u)$ in $F_{T}$.\\
It rests to verify by compactness that $(q,u)$ is a solution of the system $(SW)$.
{\hfill $\Box$}


\begin{thebibliography}{}
\bibitem{5BC}
H. Bahouri and J.-Y. Chemin, \'Equations d'ondes quasilin\'eaires et
estimation de Strichartz, Amer. J. Mathematics 121 (1999) 1337-1377.
\bibitem{5BJM}
J.-M. Bony, Calcul symbolique et propagation des singularit\'es pour
les \'equations aux d\'eriv\'ees partielles non lin\'eaires, Annales
Scientifiques de l'\'ecole Normale Sup\'erieure 14 (1981)
209-246.
\bibitem{5BG}
G. Bourdaud, R\'ealisations des espaces de Besov homog\`enes, Arkiv
fur Mathematik 26 (1998) 41-54.
\bibitem{5BD}
D. Bresch and B. Desjardins, Existence of global weak solutions for
a 2D Viscous shallow water equations and convergence to the
quasi-geostrophic model. Comm. Math. Phys., 238(1-2): 211-223,
2003.
\bibitem{5BD1}
D. Bresch and B. Desjardins, Existence of global weak solutions to
the Navier-Stokes equations for viscous compressible and heat
conducting fluids, to appear.
\bibitem{5BD2}
D. Bresch and B. Desjardins, Some diffusive capillary models of
Koretweg type. C. R. Math. Acad. Sci. Paris, Section M\'ecanique,
332(11):881-886, 2004.
\bibitem{5BDL}
D. Bresch, B. Desjardins and C.-K. Lin, On some compressible fluid
models: Korteweg, lubrification and shallow water systems. Comm.
Partial Differential Equations, 28(3-4) : 843-868, 2003.
\bibitem{5CT}
J.-Y. Chemin, Th\'eor\`emes d'unicit\'e pour le syst\`eme de
Navier-Stokes tridimensionnel, J.d'Analyse Math. 77 (1999) 27-50.
\bibitem{5CA}
J.-Y. Chemin, About Navier-Stokes system, Pr\'epublication du
Laboratoire d'Analyse Num\'erique de Paris 6 R96023 (1996).
\bibitem{5CL}
J.-Y. Chemin and N. Lerner, Flot de champs de vecteurs non
lipschitziens et \'equations de Navier-Stokes, J.Differential
Equations 121 (1992) 314-328.
\bibitem{5CK1}
H. J. Choe and H. Kim, Strong solution of the Navier-Stokes
equations for isentropic compressible fluids, J. Differential
Equations 190 (2003), 504-523.
\bibitem{5CK2}
H. J. Choe and H. Kim, Strong solution of the Navier-Stokes
equations for nonhomogeneous incompressible fluids, Math. Meth.
Appl. Sci. 28 (2005), 1-28.
\bibitem{5CR}
F. Coquel, D. Diehl, C. Merkle and C. Rohde, Sharp and diffuse
interface methods for phase transition problems in liquid-vapour
flows.
\bibitem{5DFourier}
R. Danchin, Fourier analysis method for PDE's, Preprint Novembre 2005.
\bibitem{5DG}
R. Danchin, Global Existence in Critical Spaces for Flows of
Compressible Viscous and Heat-Conductive Gases, Arch.Rational
Mech.Anal.160 (2001) 1-39
\bibitem{5DD}
Danchin.R, Density-dependent incompressible viscous fluids in critical spaces. Proc. Roy. Soc. Edimburgh Sect. A, 133(6):
1311-1334, 2003.
\bibitem{5DH}
Danchin.R, A few remarks on the Camassa-Holm equation, Differential and Integral Equations 14 (2001), 953-988.
\bibitem{5DL}
Danchin.R, Local Theory in critical Spaces for Compressible Viscous
and Heat-Conductive Gases,Communication in
Partial Differential Equations 26 (78),1183-1233 (2001)
\bibitem{5DD}
R. Danchin and B. Desjardins, Existence of solutions for
compressible fluid models of Korteweg type, Annales de l'IHP,Analyse
non
lin\'eaire 18,97-133 (2001)
\bibitem{5F}
E. Feireisl, Dynmamics of Viscous Compressible Fluids-Oxford Lecture
Series in Mathematics and its Applications-26.
\bibitem{5F2}
E. Feireisl. Compressible Navier-Stokes equations with a
non-monotone pressure law. J. Differential Equations,
184(1):97-108, 2002.
\bibitem{5F2}
E. Feireisl. On the motion of a viscous, compressible, and heat
conducting equation. Indiana Univ. Math. J.,
53(6):1705-1738, 2004.
\bibitem{5F3}
E. Feireisl, Anton\'in Novotn\'y, and Hana Petzeltov\'a. On the
existence of globally defined weak solutions to
the Navier-Stokes equations. J. Math. Fluid Mech., 3(4):358-392, 2001.
\bibitem{5H1}
David Hoff. Global existence for 1D, compressible, isentropic
Navier-Stokes equations with large initial data.
Trans. Amer. Math. Soc, 303(1):169-181, 1987.
\bibitem{5H4}
David Hoff. Discontinuous solutions of the Navier-Stokes equations
for multidimensional flows of the heat conducting fluids.
Arch. Rational Mech. Anal., 139, (1997), P. 303-354.
\bibitem{5H2}
David Hoff. Global solutions of the Navier-Stokes equations for
multidimensional compressible flow with discontinuous initial
data. J. Differential Equations, 120(1):215-254, 1995.
\bibitem{5H3}
David Hoff. Strong convergence to global solutions for
multidimensional flows of compressible, viscous fluids with
polytropic equations of state and discontinuous initial data.
Arch. Rational Mech. Anal., 132(1):1-14, 1995.
\bibitem{5HZ}
David Hoff and Kevin Zumbrum. Multi-dimensional diffusion waves for
the Navier-Stokes equations of compressible flow,
Indiana University Mathematics Journal, 1995, 44, 603-676.
\bibitem{5J}
Song Jiang and Ping Zhang. Axisymetrics solutions of the 3D
Navier-Stokes equations for compressible isentropic fluids.
J. Math. Pures Appl. (9), 82(8):949-973, 2003.
\bibitem{5K1}
A. V. Kazhikov. The equation of potential flows of a compressible
viscous fluid for small Reynolds numbers: existence,
uniqueness and stabilization of solutions. Sibirsk. Mat. Zh., 34 (1993), no. 3, p. 70-80.
\bibitem{5K}
A. V. Kazhikov and V. V. Shelukhin. Unique global solution with
respect to time of initial-boundary value problems for one-
dimensional equations of a viscous gas. Prikl. Mat. Meh., 41(2):282-291, 1977.
\bibitem{5L2}
P.-L. Lions, Mathematical Topics in Fluid Mechanics, Vol 2,
Compressible models, Oxford University Press (1998)
\bibitem{5MN}
Akitaka Matsumura and Takaaki Nishida. The initial value problem for
the equations of motion of compressible
viscous and heat-conductive fluids. Proc. Japan Acad. Ser. A Math. Sci, 55(9):337-342, 1979.
\bibitem{5MV}
A.Mellet and A.Vasseur, On the isentropic compressible Navier-Stokes
equation, Arxiv preprint math.AP/0511210, 2005 - arxiv.org
\bibitem{5Na}
J. Nash, Le probl\`eme de Cauchy pour les \'equations
diff\'erentielles d'un fluide g\'en\'eral, Bulletin de
la Soci\'et\'e Math\'ematique de France, 1962, 90, 487-497.
\bibitem{5Ro}
C. Rohde, On local and non-local Navier-Stokes-Korteweg systems for
liquid-vapour phase transitions. Technical report, Math. Institut,
Albert-Ludwigs-Universit\"at Freiburg
, 2004. Preprint.
\bibitem{5RS}
T. Runst and W. Sickel: Sobolev spaces of fractional order, Nemytskij operators, and nonlinear partial differential equations. de Gruyter Series in
Nonlinear Analysis and Applications, 3. Walter de Gruyter and Co., Berlin (1996)
\bibitem{5S}
Denis Serre. Solutions faibles globales des \'equations de
Navier-Stokes pour un fluide compressible.,
303(13):639-642, 1986
\bibitem{5So}
V.A. Solonnikov. Estimates for solutions of nonstationary
Navier-Stokes systems. Zap. Nauchn. Sem. LOMI, 38,
(1973), p.153-231; J. Soviet Math. 8, (1977), p. 467-529.
\bibitem{5V}
V. Valli, W. Zajaczkowski. Navier-Stokes equations for compressible
fluids: global existence and qualitative properties of the solutions
in the general case. Commun. Math. Phys., 103 (1986) no 2., p.
259-296.
\bibitem{5W}
Weike Wang and Chao-Jiang Xu. The Cauchy problem for viscous shallow
water equations.
Rev. Mat. Iberoamericana  21, no. 1 (2005), 1-24.
\end{thebibliography}
\end{document}